\documentclass{amsart}
\usepackage[english]{babel}
\usepackage{amsaddr}

\usepackage[utf8]{inputenc}
\usepackage[T1]{fontenc}
\usepackage{csquotes}
\usepackage{enumitem}
\usepackage{esint}

\usepackage[backend=biber]{biblatex}
\usepackage{lmodern}

\usepackage{amsmath}

\usepackage{amssymb}

\usepackage{bbm}

\usepackage{fancyhdr}

\numberwithin{equation}{section}

\usepackage{xcolor}
\usepackage{geometry}
\newtheoremstyle{theorem2}
{8pt}
{8pt}
{\itshape}
{}
{\bfseries}
{.}
{.5em}
{}

\newtheoremstyle{definition2}
{8pt}
{8pt}
{\itshape}
{}
{\bfseries}
{.}
{.5em}
{}

\theoremstyle{theorem2}

\newtheorem*{theorem*}{Theorem}

\newtheorem{lemma}{Lemma}[section]
\newtheorem{teo}[lemma]{Theorem}

\newtheorem{cor}[lemma]{Corollary}

\theoremstyle{definition2}

\newtheorem{deff}[lemma]{Definition}

\theoremstyle{plain}
\newtheorem*{mainthm}{Main Theorem}

\newcommand{\orcid}[1]{ORCID: #1}

\newcommand{\Dfdelta}[1]{{D}_{#1}f_\delta}
\newcommand{\Dfdeltaij}{D^2_{i,j}f_\delta}
\newcommand{\Dfdeltajnl}{D^2_{j,l+n}f_\delta}
\newcommand{\Dfdeltajln}{D^2_{j,l-n}f_\delta}

\begin{document}


\title{Lipschitz regularity for solutions to an orthotropic $q$-Laplacian-type equation in the Heisenberg group}

\author{Michele Circelli (\orcid{0009-0005-4928-8310})}
\address{University of Bologna, Department of Mathematics, 40126-Bologna (Italy)}
\email{michele.circelli2@unibo.it}

\author{Giovanna Citti (\orcid{0000-0003-2316-3315})}
\address{University of Bologna, Department of Mathematics, 40126-Bologna (Italy)}
\email{giovanna.citti@unibo.it}

\author{Albert Clop (\orcid{0000-0002-0187-6288})}
\address{University of Barcelona, Department of Mathematics and Computer Science, 08007-Barcelona (Catalonia)}
\email{albert.clop@ub.edu}

\begin{abstract}
We establish the local Lipschitz regularity for solutions to an orthotropic $q$-Laplacian-type equation within the Heisenberg group. Our approach is largely inspired by the works of X. Zhong, who investigated the $q$-Laplacian in the same setting and proved the Hölder regularity for the gradient of solutions. Due to the degeneracy of the current equation, such regularity for the gradient of solutions is not even known in the Euclidean setting for dimensions greater than 2, where only boundedness is expected.
\end{abstract}

\maketitle

\section{Introduction}

The standard $q$-Laplacian equation in the Euclidean setting $\mathbb{R}^n$ is a second-order, divergence-type equation, expressed as 
\begin{equation}\label{divergence}
\textnormal{div }( D f(D u) )=0,
\end{equation}
where $f:\mathbb{R}^n\to\mathbb{R}$ is given by $f(z)= \frac{1}{q}|z|^{q}$ and $D$ denotes the Euclidean gradient. It is well known that the optimal regularity for solutions is the Hölder continuity of the gradient. This classical result was obtained by N. Uraltseva \cite{Uraltseva} and K. Uhlenbeck \cite{Uhlenbeck} for $q>2$, and by E. DiBenedetto \cite{DiBenedetto} and P. Tolksdorf \cite{Tolksdorf} for general values of $q$.

A degenerate version of this equation was introduced in \cite{Uralurdal} by choosing 
\begin{equation}\label{ortotropic1}
f(z)=\sum_{i=1}^{n}\frac{1}{q_i}|z_i|^{q_i},\quad\textnormal{with }1<q_1\leq\ldots\leq q_n,
\end{equation} 
in equation \eqref{divergence}. The authors provided a proof of the local boundedness for the gradient of \textit{locally bounded} solutions when $q_1\geq4$ and $\frac{q_n}{q_1}<2$. This result was improved, with no limitations on the ratio $\frac{q_n}{q_1}$, in \cite{bousquet2020lipschitz} for $q_1\geq 2$ and in \cite{bousquet2023singular} for $1<q_1\leq\ldots\leq q_n\leq2$.

If $q_1=\ldots =q_n$, then
\begin{equation}\label{ortotropic}
f(z)=\frac{1}{q}\sum_{i=1}^{n}|z_i|^{q},
\end{equation}
and in this case, \eqref{divergence} has received significant attention since the paper \cite{brasco2013congested}, where the authors related the equation to the Beckmann formulation of optimal transportation with congestion effects. In \cite{bousquet2018lipschitz}, the aforementioned regularity results for solutions were extended to a more degenerate equation for $q\geq2$, without any boundedness assumption on the solutions themselves. See also \cite{bousquet2014lipschitz} for intermediate results and \cite{Demengel} for an alternative proof based on viscosity methods.

In the particular case $n=2$, Bousquet and Brasco \cite{bousquet2018c1} proved that weak solutions of \eqref{divergence} are $C^1$ for $1<p<\infty$. Moreover, derivatives of solutions have a logarithmic modulus of continuity: see \cite{ricciotti2018regularity} for the case $1<q<2$, and \cite{lindqvist2018regularity} for the case $q\geq2$ (the latter result holds for the more general equation \eqref{ortotropic}). See also \cite{BrascoCarlier2014}, \cite{bousquet2014lipschitz}, and \cite{brasco2017sobolev} for intermediate results.

However, the problem of regularity remains largely open when $n>2$. The equation is more degenerate than the $q$-Laplacian since the principal part vanishes on a larger set. Consequently, it is not yet clear whether the gradient is merely bounded or H\"older continuous, even in the Euclidean setting.

\medskip
The Heisenberg group $\mathbb{H}^n$ is a connected and simply connected Lie group, whose tangent bundle can be represented as $H\mathbb{H}^n \oplus V\mathbb{H}^n$, called the \textit{horizontal} and vertical bundles, respectively. A choice of a metric on $H\mathbb{H}^n$ allows for the selection of an orthonormal frame $X_1,\ldots,X_{2n}$ for $H\mathbb{H}^n$. Additionally, it allows for the definition of a sub-Riemannian distance on the space: its associated balls will be denoted by $B_{SR}$. We will denote the horizontal gradient of a scalar function $u$ by $\nabla_Hu=\left(X_1u,\ldots,X_{2n}u\right)$, and the horizontal divergence of a vector-valued function $\phi = \sum_{i=1}^{2n} \phi_i X_i$ by $\textnormal{div}_H \phi = \sum_{i=1}^{2n} X_i \phi_i$. Consequently, if $f:\mathbb{R}^{2n}\to\mathbb{R}$, then equation \eqref{divergence} becomes in this setting
\begin{equation}
\label{Hdivergence}
\textnormal{div}_H ( D f(\nabla_H u) )=0,
\end{equation}
where $Df$ denotes the Euclidean gradient of $f$ in $\mathbb{R}^{2n}$, and the natural domain for its solutions is the Sobolev space $HW^{1,q}(\Omega)$ associated with $\nabla_H$.

As before, if $f(z)=\frac{1}{q}|z|^q$, then \eqref{Hdivergence} becomes the standard $q$-Laplace equation in $\mathbb{H}^n$. The optimal regularity for solutions is the Hölder continuity of the horizontal derivatives, established in \cite{zhong2017regularity} for $2\leq q<\infty$, and in \cite{mukherjee20211} for $1<q<2$. See also \cite{ricciotti2018} for an alternative proof in $\mathbb{H}^1$ for $q>4$, \cite{cittimukherjee} for a generalization of the equation, and \cite{ricciotti} for a general overview on this topic.

Recently, the congested optimal transport problem has been introduced in the Heisenberg group, thus leading to degenerate versions of the $q$-Laplacian equation in this setting (see \cite{circelli2024continuous} and \cite{circelli2023transport}). The Heisenberg group can also be identified as a real hypersurface in $\mathbb{C}^{n+1}$. With this interpretation, the horizontal tangent space at any point is invariant under multiplication by the imaginary unit, which couples $X_i$ and $X_{i+n}$. Correspondingly, the analogous function $f:\mathbb{R}^{2n}\to\mathbb{R}$ defined in \eqref{ortotropic} is
\begin{equation}\label{ortotropicH}
f(z)=\frac{1}{q}\sum_{i=1}^{n}\left(z_i^2+z_{i+n}^2\right)^{\frac{q}{2}}.
\end{equation}
Denoting by 
\begin{equation*}
\aligned
    &\lambda_i(z)=\left(z_i^2+z_{i+n}^2\right)^{\frac{q-2}{2}},\quad \forall i\in\{1,\ldots,n\},\\
    &\lambda_{i}(z) = \left(z_{i-n}^2+z_i^2\right)^{\frac{q-2}{2}},\quad \forall i\in\{n+1,\ldots,2n\},
\endaligned
\end{equation*}
the associated equation has the following orthotropic $q$-Laplacian-type structure: 
\begin{equation}
\label{eq}
\textnormal{div}_H ( D f(\nabla_H u) ) = \sum_{i=1}^{2n}X_i(\lambda_i(\nabla_H u) X_iu)=0,
\end{equation}
where $Df=(D_1f, \dots, D_{2n}f)$ is the Euclidean $2n$-dimensional gradient of $f$.

The $q$-Laplace equation degenerates at $z=0$, while equation \eqref{eq} is much more degenerate, as it degenerates in the unbounded set
\begin{equation*}
	\bigcup_{i=1}^n\{ \lambda_i (z) =0\} =\bigcup_{i=1}^n\left\{z_i^2+z_{i+n}^2=0\right\}\subset\mathbb{R}^{2n},
\end{equation*}
which is the union of $2n-2$ dimensional submanifolds. For this reason, we do not expect the Hölder regularity of the gradient of solutions, as in the standard $q$-Laplacian case \cite{zhong2017regularity}, but only the boundedness of the gradient, in analogy with the Euclidean setting \cite{bousquet2018lipschitz} or \cite{Demengel}.

In this paper we address this regularity result in the range $2 \le q < +\infty$, by adapting techniques introduced by Zhong for the standard (non-orthotropic) $q$-Laplacian in $\mathbb{H}^n$; see \cite{zhong2017regularity}. In that setting the method applies to the full range $1<q<+\infty$, and in particular yields Hölder regularity of the horizontal gradient for $q \ge 2$. In the orthotropic case considered here, however, the adaptation of strategy \cite{zhong2017regularity} uses in a crucial way the assumption $q \ge 2$, and our arguments do not extend to the singular regime $1<q<2$. As a consequence, even the Lipschitz regularity of weak solutions to our orthotropic equation remains, to the best of our knowledge, an open problem when $1<q<2$. We also point out that, for the standard $q$-Laplacian, the analysis of Hölder regularity of the horizontal gradient of solutions in the singular range $1<q<2$ requires substantially different techniques, see \cite{mukherjee20211}, and it is therefore natural to expect that the corresponding orthotropic problem will require new ideas, rather than a straightforward modification of the arguments developed in the present paper. A systematic study of the regularity theory in the range $1<q<2$ is therefore left for future work.

The main theorem of this paper is the following one.
\begin{mainthm}
Let $2\leq q<\infty$ and $u\in HW^{1,q}(\Omega)$ be a weak solution of \eqref{eq}. Then, $\nabla_Hu\in L^\infty_{loc}(\Omega,\mathbb{R}^{2n})$. Moreover, for any ball $B_{SR}(r)$ such that $B_{SR}(2r)\subset\Omega$ it holds that 
\begin{equation*}
	\|\nabla_H u\|_{L^\infty(B_{SR}(r))}\leq c \left(\fint_{B_{SR}(2r)}|\nabla_H u|_H^qdx\right)^{\frac{1}{q}},
\end{equation*}
where $c=c(q,N,L)>0$ and $N=2n+2$ is the homogeneous dimension of $\mathbb{H}^n$. 
\end{mainthm}

For analogous results in the parabolic setting, in the range $2\leq q\leq 4$, see \cite{circelli2025gradientestimatesorthotropicnonlinear}. 

Let us finally point out that when $q=2$, equation \eqref{eq} becomes the standard sub-Laplacian $\Delta_Hu:=\sum_{i=1}^{2n}X_i^2u=0$, for which the smoothness of solutions directly follows from \cite{Hormander1967}. 

\medskip
The main difficulty in studying regularity for PDEs in this setting arises from the non commutativity of derivatives: $[X_i,X_{i+n}]=X_{2n+1}$. Note that $X_i$ and $X_{i+n}$ are elements of the horizontal gradient and therefore first derivatives, while their commutator $X_{2n+1}$ acts as a second-degree derivative. Due to this non-commutativity, the equation satisfied by the first derivatives (necessary to estimate the gradient) involves a derivative of solutions in the direction of $X_{2n+1}$, leading to an higher-order term.

As is common in the study of regularity theory for PDEs, we first approximate the equation to make it non-degenerate. Following the approach in \cite{capogna2020regularity}, \cite{capogna2021lipschitz}, \cite{capogna2023regularity}, and \cite{capogna2019conformality}, we introduce a Riemannian approximation of the non-degenerate equation to achieve smoothness of solutions. We then obtain estimates for derivatives of approximating solutions, independent of the regularizing parameters: in the spirit of  \cite{zhong2017regularity}, we introduce mixed Caccioppoli-type inequalities for the gradient of approximating solutions, which involve both horizontal and vertical derivatives. Interestingly, this same technique is a key tool used in \cite{bousquet2018lipschitz} in the Euclidean setting. By combining ideas from both papers, we derive a standard Caccioppoli-type inequality for the gradient of approximating solutions. Using the Moser iteration scheme, we obtain a uniform bound (with respect to the approximating parameters) for the $L^\infty$ norm of the gradient of solutions to the approximating equation. Passing to the limit, we establish a bound on the $L^\infty$ norm for the horizontal gradient of solutions to the original equation, thereby proving local Lipschitz continuity.

We finally point out that our arguments make essential use of the step-two structure of $\mathbb{H}^n$ and of the single non-trivial commutation relation $[X_i,X_{i+n}]=X_{2n+1}, \qquad i=1,\dots,n.$ This specific structure plays a crucial role in the derivation of the mixed Caccioppoli-type estimates used to control the vertical derivative $X_{2n+1}u$ in terms of the horizontal derivatives. This suggests a possible extension of our results to more general step-two Carnot groups. However, for such groups the commutator structure is more involved, and an extension of our method would require a finer analysis.

The structure of the paper is as follows. In Section \ref{preliminaries}, we recall some preliminary definitions and results, introduce the Heisenberg group, and review how it arises as limit of Riemannian manifolds. In Section \ref{Riemannianapproxschemesection}, we construct the (Riemannian) uniformly elliptic equation that approximates the original one. Section \ref{Lipschitzsection} is devoted to the study of the local Lipschitz continuity of solutions.

\section{Preliminaries}\label{preliminaries}
\subsection{The Heisenberg group $\mathbb{H}^n$}

Let $n\geq1$. The $n$-th \textit{Heisenberg group} $\mathbb{H}^n$ is the connected and simply connected Lie group, whose Lie algebra $\mathfrak{h}^n$ is stratified of step 2: i.e. it is the direct sum of two linear subspaces
$$ \mathfrak{h}^n= \mathfrak{h}^n_1 \oplus \mathfrak{h}^n_2,$$
where $\dim(\mathfrak{h}^n_1)= 2n$, 
$\dim(\mathfrak{h}^n_2)= 1,$ $[\mathfrak{h}^n_1, \mathfrak{h}^n_1] = \mathfrak{h}^n_2$ and 
$[\mathfrak{h}^n_1, \mathfrak{h}^n_2] = 0.$
Due to the stratification, we will assign degree 1 to the elements of $\mathfrak{h}^n_1$,  degree 2 to $\mathfrak{h}^n_2$, so that the  homogeneous dimension of $\mathbb{H}^n$ is 
\begin{equation*}\label{homogdim}
N:=\sum_{i=1}^2i\dim(\mathfrak{h}^n_i)=2n+2.
\end{equation*}
Since the Lie algebra is nilpotent, it is well-known that the exponential map $\mathrm{exp}: \mathfrak{h}^n\to \mathbb{H}^n$ is a  global diffeomorphism. Hence, if we choose a  basis $X_1,\ldots,X_{2n}$ of $\mathfrak{h}^n_1$, and a vector $X_{2n+1}\in \mathfrak{h}^n_2$, such that $[X_i,X_{i+n}]=X_{2n+1}$, we can to introduce on $\mathbb{H}^n$ suitable choice of exponential coordinates
\begin{equation*}
\mathbb{H}^n\ni\exp\left(x_1X_1+x_2X_2+\ldots+x_{2n}X_{2n}+x_{2n+1}X_{2n+1}\right)\longleftrightarrow\left(x_1,\ldots,x_{2n+1}\right)\in\mathbb{R}^{2n+1}.
\end{equation*} 
In this system of coordinates the group law will be given, through the Baker-Campbell-Hausdorff formula, by  
\begin{equation*}\label{grouplaw}
	x \cdot y := \bigg{(}x_1+y_1,\ldots,x_{2n}+y_{2n}, x_{2n+1}+y_{2n+1}+\frac{1}{2}\sum_{i=1}^{2n} (x_iy_{i+n}-x_{i+n}y_i)\bigg{)},
\end{equation*}
the unit element will be denoted by $0$. In addition, the vector fields $X_1,\ldots,X_{2n+1}$, left invariant w.r.t. \eqref{grouplaw}, in coordinates read as
\begin{equation}
	X_i := \partial_{i} -\frac{x_{i+n}}{2} \partial_{2n+1} \,, i=1,\dots,n,\quad 
		X_{i+n} := \partial_{i+n}+\frac{x_i}{2} \partial_{2n+1},\ \ i=1,\dots,n,
\end{equation}
See Section 2 in \cite[Section 2, Theorem 2.2.18]{Bonfiglioli} for more details about this topic.

The horizontal layer $\mathfrak{h}^n_1$ of the algebra defines a sub-bundle $H\mathbb{H}^n$ of the tangent bundle $T\mathbb{H}^n$, whose fibre at a point $x \in \mathbb{H}^n$ is
$$H_x \mathbb{H}^n=\mathrm{span} \{ X_1(x), \dots, X_n(x), X_{n+1}(x), \dots, X_{2n}(x) \}.$$

One can equip $\mathbb{H}^n$ with a left-invariant sub-Riemannian metric as follows: we fix an inner product $\langle \cdot, \cdot \rangle_H$ on $\mathfrak{h}_1^n$ that makes $\{X_1, \dots, X_n, X_{n+1}, \dots, X_{2n}\}$ an orthonormal basis and we denote by $|\cdot|_H$ the norm associated with such a scalar product. We keep the same notation both for the corresponding scalar product and norm on each fibre $H_x\mathbb{H}^n$, $x\in\mathbb{H}^n$.

We call horizontal curve $\sigma$  any absolutely continuous curve $\sigma\in AC([0,1],\mathbb{R}^{2n+1})$, whose tangent vector $\dot{\sigma}(t)$ belongs to $H_{\sigma(t)}\mathbb{H}^n$ at almost every $t\in[0,1]$. 
Due to the stratification of the space the H\"ormander condition is satisfied, and the  Rashevsky-Chow's theorem guarantees that any couple of points can be joined with an horizontal curve  \cite{Chow}. It is then possible to give the following definition of distance. Given $x,y\in\mathbb{H}^n$, the Carnot-Carathéodory distance between them is defined as
\begin{equation*}\label{CCdistance}
	d_{SR}(x,y):= \inf \  \left\{  \int_0^1|\dot{\sigma}(t)|_H dt \ | \  \sigma \text{ is horizontal}, \ \sigma(0)=x, \ \sigma(1)=y \right\}. 
\end{equation*}

If $x\in\mathbb{H}^n$ and $r>0$, the corresponding ball will be denoted  by
\begin{equation*}
	B_{SR}(x,r):=\left\{y\in\mathbb{H}^n: d_{SR}(x,y)\leq r\right\}.
\end{equation*}

The Lebesgue measure $\mathcal{L}^{2n+1}$ is a Haar measure on $\mathbb{H}^n\simeq\mathbb{R}^{2n+1}$; we will often denote it by $dx$.

For any open set $\Omega\subseteq\mathbb{H}^n$ and any smooth enough function $u:\Omega\to\mathbb{R}$, we denote by 
$$\nabla_Hu=\sum_{i=1}^{2n}X_iuX_i$$ 
its horizontal gradient. 
For every $1\leq q\leq\infty$, the space
\begin{equation}\label{horsob}
	HW^{1,q}(\Omega):=\left\{u\in L^{q}(\Omega):\nabla_Hu\in L^{q}(\Omega, H\Omega)\right\},
\end{equation}
is a Banach space equipped with the norm
\begin{equation*}
	\|u\|_{HW^{1,q}(\Omega)}:=\|u\|_{L^q(\Omega)}+\|\nabla_H u\|_{L^q(\Omega,H\Omega)}.
\end{equation*}
Moreover, for $1\leq q<\infty$, we denote by
\begin{equation*}
	HW_0^{1,q}(\Omega):=\overline{C_0^{\infty}(\Omega)}^{HW^{1,q}(\Omega)}.
\end{equation*}

\subsection{Riemannian approximation of the sub-Riemannian Heisenberg group}

The Heisenberg group $\left(\mathbb{H}^n,d_{SR}\right)$ arises as the pointed Hausdorff-Gromov limit of Riemannian manifolds, in which the non-horizontal direction is increasingly penalized.
Precisely, for any $\epsilon>0$ one can consider the metric tensor $g_\epsilon$ that makes 
\begin{equation*}
	X_1(x),\ldots,X_{2n}(x),\epsilon X_{2n+1}(x)
\end{equation*}
an orthonormal basis for $T_x\mathbb{H}^n$, at any point $x\in \mathbb{H}^n$. We relabel 
\begin{equation*}
	X_i^\epsilon:=X_i,\forall i=1,\ldots,2n\text{ and }X_{2n+1}^\epsilon:=\epsilon X_{2n+1},
\end{equation*}
and we denote by $|\cdot|_\epsilon:=\sqrt{g_\epsilon(\cdot,\cdot)}$ and by $d_\epsilon$ the associated control distance: it turns out to be left-invariant with respect to \eqref{grouplaw}, since $X^\epsilon_1,\ldots, X^\epsilon_{2n}, X_{2n+1}^\epsilon$ are themselves left invariant. It has been proved by Gromov in \cite{Gromov} that for  any couple of points $x,y\in \mathbb{H}^n$ the distance satisfies
\begin{equation*}
	d_{SR}(x,y)=\lim_{\epsilon\to0}d_\epsilon(x,y),
\end{equation*}
and the convergence is uniform on compact subsets of $\mathbb{H}^n\times\mathbb{H}^n$, see also \cite{capogna2016regularity}. In particular the Riemannian balls $B_\epsilon$, associated with $d_\epsilon$, satisfy $B_\epsilon\to B_{SR}$, in the Hausdorff-Gromov sense. Moreover, one can consider a regularized gauge function
\begin{equation*}
	G_\epsilon^2(x):=\sum_{i=1} ^{2n} |x_i|^2+ \min\left\{\frac{|x_{2n+1}|^2}{\epsilon^2}, |x_{2n+1}|\right\}, \quad x\in \mathbb{H}^n,
\end{equation*}
and \cite[Lemma 2.13]{capogna2016regularity} implies that the distance function $d_{G,\epsilon}(x,y):=G_\epsilon(y^{-1}\cdot x)$ is equivalent to $d_\epsilon$, i.e. there exists a constant $A>0$ such that $\forall \epsilon>0$ and $\forall x,y\in \mathbb{H}^n$ 
\begin{equation}\label{gauge}
	A^{-1}d_{G,\epsilon}(x,y)\leq d_\epsilon (x,y)\leq Ad_{G,\epsilon}(x,y).
\end{equation}
In particular, \eqref{gauge} implies that the doubling property for Riemannian balls
\begin{equation*}
	\mathcal{L}^{2n+1}(B_\epsilon(x,2r))\leq C\mathcal{L}^{2n+1}(B_\epsilon(x,r)),
\end{equation*}
holds uniformly in $\epsilon$, with a constant $C>1$ independent of $\epsilon$. 

If $\Omega\subseteq \mathbb{H}^n$ is an open set and $u:\mathbb{H}^n\to\mathbb{R}$ is any smooth enough function, then the gradient associated with the Riemannian metric $g_\epsilon$ is
\begin{equation*}
	\nabla_\epsilon u:=\sum_{i=1}^{2n+1}X_i^\epsilon uX_i^\epsilon=\sum_{i=1}^{2n}X_iuX_i+\epsilon^2X_{2n+1}uX_{2n+1}.
\end{equation*}

For a regular vector field $\phi= \sum_{i=1}^{2n+1} \phi_i X_i^\epsilon$ will also denote 
$$\textnormal{div}_\epsilon\phi = \sum_{i=1}^{2n+1} X_i^\epsilon \phi_i.$$
Formally we have
\begin{equation*}
	\nabla_\epsilon u\longrightarrow \nabla_Hu,\quad \textnormal{div}_\epsilon\phi \longrightarrow \textnormal{div}_H \phi,\quad\text{ as }\epsilon\to0.
\end{equation*}

For every $1\leq q\leq\infty$, one can define the Banach space
\begin{equation}\label{riemanniansob}
	W_\epsilon^{1,q}(\Omega):=\left\{u\in L^{q}(\Omega):\nabla_\epsilon u\in L^{q}(\Omega,T\Omega)\right\},
\end{equation}
equipped with the norm
\begin{equation*}
	\|u\|_{W_\epsilon^{1,q}(\Omega)}:=\|u\|_{L^q(\Omega)}+\|\nabla_\epsilon u\|_{L^q(\Omega,T\Omega)}.
\end{equation*}

Finally, let us recall that Sobolev's inequalities hold with constants independent of $\epsilon$, see for instance \cite{capogna2016regularity}.

\begin{teo}\label{sobolevineqriemannian}
Let $1\leq p<N$. For any Riemannian ball $B_\epsilon(r)\subset\mathbb{H}^n$ and any $u\in W^{1,p}_\epsilon(B_\epsilon(r))$, it holds
\begin{equation*}
	\left(\fint_{B_\epsilon(r)}|u|^{\frac{Np}{N-p}}dx\right)^{\frac{N-p}{Np}}\leq cr\left(\fint_{B_\epsilon(r)}|\nabla_\epsilon u|_\epsilon^pdx\right)^{\frac{1}{p}},
\end{equation*}
where $c=c(p,N)$ is independent of $\epsilon$.
\end{teo}

\section{The orthotropic $q$-Laplacian-type equation}\label{Riemannianapproxschemesection}
Let $q\geq2$ and let us consider the equation
\eqref{eq}. 
A direct computation shows that  the following growth condition holds
\begin{equation*}\label{structurecondition1}
    \sum_{i=1}^{n} \lambda_i(z)(\xi_i^2+\xi_{i+n}^2)\leq\langle D^2f(z)\xi,\xi\rangle\leq (q-1)\sum_{i=1}^{n}\lambda_i(z)(\xi_i^2+\xi_{i+n}^2),
\end{equation*}
for any $\xi\in\mathbb{R}^{2n}$.

\begin{deff}[Weak Solution]
	We say that a function $u\in HW^{1,q}(\Omega)$ is a weak solution to the equation \eqref{eq} if 
	\begin{equation*}
		\sum_{i=1}^{2n}\int_{\Omega}D_if(\nabla_H u)X_i\varphi dx=0,
	\end{equation*}
	for any $\varphi\in C^\infty_0(\Omega)$.
\end{deff}

\subsection{Riemannian approximation scheme}
The approximation procedure we introduce below is a regularization scheme widely used in the literature to prove a-priori estimates for weak solutions of subelliptic PDEs, see for instance \cite{capogna2020regularity}, \cite{capogna2021lipschitz} and \cite{capogna2023regularity}.
Then, for any $\delta\in(0,1)$ we denote by $f_\delta:\mathbb{R}^{2n+1}\to\mathbb{R}$ the function
$$f_{\delta}(z) = \frac{1}{q}\sum_{i=1}^{n} \left(\delta+z_i^2+z_{i+n}^2+ z_{2n+1}^2\right)^{\frac{q}{2}} $$
and we consider the differential equation 
\begin{equation}\label{eqneriemannianapprox}
\textnormal{div}_\epsilon ( D f_\delta(\nabla_\epsilon u) )=0,
\end{equation}
where $Df_\delta$ denotes the Euclidean gradient of $f_\delta$ in $\mathbb{R}^{2n+1}$.

If we denote by 
\begin{equation}\label{lambdaideltaeps}
\begin{split}
    &\lambda_{i,\delta}(z):=(\delta+z_i^2 + z_{i+n}^2+z_{2n+1}^2)^{(q-2)/2},\quad \forall i\in\{1,\ldots,n\},\\
    &\lambda_{i,\delta}(z):=(\delta+z_{i-n}^2 + z_{i}^2+z_{2n+1}^2)^{(q-2)/2},\quad \forall i\in\{n+1,\ldots,2n\},\\&\lambda_{2n+1, \delta}(z):= \sum_{i=1}^{n}\lambda_{i, \delta}(z),
\end{split}
\end{equation}
then the gradient of $f_\delta$ can be written as follows
\begin{equation}\label{structurecondition02}
    Df_{\delta}(z) = \Big(\lambda_{1,\delta}(z) z_1, \cdots, \lambda_{2n,\delta}(z) z_{2n},  \lambda_{2n+1,\delta}(z) z_{2n+1}\Big).
\end{equation}
Hence, the explicit expression of the regularized equation 
becomes 
$$\sum_{i=1}^{2n+1}X_i^\epsilon (\lambda_{i, \delta}(\nabla_\epsilon u) X_i^\epsilon u)=0. 
$$

\begin{deff}[$\epsilon$-Weak Solution]
	We say that a function $u_\epsilon\in W_\epsilon^{1,q}(\Omega)$ is a weak solution to the equation \eqref{eqneriemannianapprox} if 
	\begin{equation}\label{weakformulation}
		\sum_{i=1}^{2n+1}\int_{\Omega}D_i f_{\delta}(\nabla_\epsilon u_\epsilon)X_i^\epsilon\varphi dx=0,
	\end{equation}
	for any $\varphi\in C^\infty_0(\Omega)$.
\end{deff}

Finally, the structure condition becomes
\begin{equation}\label{structurecondition01}
	\sum_{i=1}^{n} \lambda_{i,\delta}(z)(\xi_i^2 + \xi_{i+n}^2 + \xi_{2n+1}^2)\leq\langle D^2f_{\delta}(z)\xi,\xi\rangle \leq L\sum_{i=1}^{n} \lambda_{i,\delta}(z)(\xi_i^2 + \xi_{i+n}^2+ \xi _{2n+1}^2),    
\end{equation}
for any $\xi\in\mathbb{R}^{2n+1}$, where $L=L(q,n)>1$ is a constant. 
Hence, from the Euclidean and Riemannian elliptic regularity theory (see for instance \cite{Uraltsevabook}) it follows that every $\epsilon$-weak solution $u_\epsilon$ to \eqref{eqneriemannianapprox} is smooth, that is $u_\epsilon\in C^\infty(\Omega)$.

\section{Local Lipschitz regularity for solutions}\label{Lipschitzsection}

In this section we prove the main theorem of this paper. The argument proceeds as follows: we first establish a Caccioppoli-type inequality for the first derivatives of solutions to the regularized equation; at this level, the estimate still depends on the vertical derivative of the solution itself. We then remove this dependence by adapting to the present orthotropic setting a mixed Caccioppoli-type inequality, originally introduced by Zhong in \cite{zhong2017regularity} for solutions to the standard $q$-Laplacian in $\mathbb{H}^n$. Finally, combining these estimates with a Moser iteration scheme yields an $L^\infty$ bound for the Riemannian gradient of approximating solutions, and passing to the limit in the regularization parameters leads to the desired local Lipschitz regularity for weak solutions of the original equation.

\subsection{Uniform Caccioppoli-type inequalities for first derivatives of $\epsilon$-weak solutions}

The aim of this subsection is to get higher regularity estimates for weak solutions $u_\epsilon$ to \eqref{eqneriemannianapprox} that are stable in $\epsilon$ and $\delta$.
Through the whole subsection, with an abuse of notation, we will drop the index $\epsilon$ and we will denote by $u$ a weak solution to \eqref{eqneriemannianapprox}. Moreover, throughout the rest of the paper we denote by $c$ a positive constant, that may vary from
line to line. Except explicitly being specified, it depends only on the dimension $n$ and on the constants $q$ and $L$ in the structure condition \eqref{structurecondition01}. However, it does not depend on the approximating parameters $\epsilon$ and $\delta$.

For every $v\in C^\infty(\Omega)$ we define the following two quantities that play the formal role of the weighted gradient, and its formal gradient square, respectively: 
    
\begin{equation}\label{26ottobre1}
    \aligned
    &\Psi(v)=\Psi_{\epsilon,\delta}(v) =\Big(\lambda_{1,\delta}(\nabla_\epsilon v) X_1^\epsilon v, \cdots, \lambda_{2n+1,\delta}(\nabla_\epsilon v) X_{2n+1}^\epsilon v\Big),\\
    &||\Psi(v)||^2_\lambda:=\sum_{i=1}^n \lambda_{i,\delta}(\nabla_\epsilon v) ((X_i^\epsilon v)^2+(X_{i+n}^\epsilon v)^2+(X_{2n+1}^\epsilon v)^2).
    \endaligned
\end{equation}

Given a weak solution $u\in C^\infty(\Omega)$ to \eqref{eqneriemannianapprox}, the standard Cacciopoli inequality for the first derivatives is an estimate of 
$$ ||\Psi(X_j^\epsilon u )||^2_\lambda :=  
\sum_{i=1}^n \lambda_{i,\delta}(\nabla_\epsilon u) ((X_i^\epsilon X_j^\epsilon u)^2+(X_{i+n}^\epsilon X_j^\epsilon u)^2+(X_{2n+1}^\epsilon X_j^\epsilon u)^2),$$ 
where $j=1,\ldots,2n+1$. The idea is to exchange the role of $X_i^\epsilon$ and $X_j^\epsilon$ and estimate instead 
$$ ||\Phi(X_j^\epsilon u )||^2_\lambda :=  
\sum_{i=1}^n \lambda_{i,\delta}(\nabla_\epsilon u) ((X_j^\epsilon X_i^\epsilon u)^2+(X_{j}^\epsilon X_{i+n}^\epsilon u)^2+(X_{j}^\epsilon X_{2n+1}^\epsilon u)^2),$$ 
so that 
\begin{equation}\label{Phi_def}
    ||\Phi(\nabla_\epsilon u)||^2_\lambda := \sum_{j=1}^{2n+1} 
||\Phi(X^\epsilon_j u )||^2_\lambda.
\end{equation}

We also denote by
\begin{equation}\label{26ottobre}
    \aligned
    \Lambda=\Lambda_{\epsilon,\delta}&:=\sum_{i=1}^n \lambda_{i,\delta}(\nabla_\epsilon u),\\
    w(\nabla_\epsilon u)&:=\left(\delta+|\nabla_\epsilon u|_\epsilon^2\right)^\frac{1}{2}.
    \endaligned
\end{equation}

\begin{lemma}\label{caccioppoliXuriemannian}[Caccioppoli for $\nabla_\epsilon u$]
There exists a constant $c=c(q,n,L)>0$, independent of $\epsilon$ and $\delta$, such that, for every weak solution $u\in W_\epsilon^{1,q}(\Omega)$ to \eqref{eqneriemannianapprox}, for every $\beta\geq0$ and for every non-negative $\eta\in C_0^\infty(\Omega)$ one has
\begin{equation*}
\aligned
    &\int_\Omega\eta^2\  w(\nabla_\epsilon u)^\beta\   ||\Phi(\nabla_\epsilon u )||^2_\lambda\   dx \\&\leq c(\beta+1)\int_\Omega (|\nabla_\epsilon\eta|_\epsilon^2+\eta|X_{2n+1}\eta|)\  w(\nabla_\epsilon u)^{\beta+2}\ \Lambda dx\\&+c(\beta+1)^2\int_\Omega\eta^2\  w(\nabla_\epsilon u)^\beta\ |X_{2n+1}u|^2\ \Lambda dx,
\endaligned
\end{equation*}
where $\|\Phi\|_\lambda^2$ is defined in \eqref{Phi_def} and $w$ and $\Lambda$ are defined in \eqref{26ottobre}.

\end{lemma}

\begin{proof}
We use $\varphi=X^\epsilon_l(\eta^2w^\beta X^\epsilon_lu)$, with $\beta\geq0$ and $l\in\{1,\ldots,2n+1\}$, as test function in the weak formulation \eqref{weakformulation} of the equation \eqref{eqneriemannianapprox}
\begin{equation}\label{weakform}
	0=\int_\Omega\sum_{i=1}^{2n+1}{\Dfdelta i} (\nabla_\epsilon u)X_i^\epsilon\varphi dx =\int_\Omega\sum_{i=1}^{2n+1}{\Dfdelta i} (\nabla_\epsilon u)X_i^\epsilon(X^\epsilon_l(\eta^2w^\beta X^\epsilon_lu))dx. \nonumber
\end{equation}
Now we assume that $l\in\{1,\ldots,n\}$ and we exchange the order of differentiation in the second term. Using the fact that $l\leq n$, so that $X_i^\epsilon X^\epsilon_l 
 =-\delta_{i,l+n} X_{2n+1} + X_l^\epsilon X^\epsilon_i$
 we get 
\begin{equation*}
	0=-\sum_{i=1}^{2n+1}\int_\Omega X_l^\epsilon({\Dfdelta i} (\nabla_\epsilon u))X_i^\epsilon(\eta^2w^\beta X^\epsilon_lu)dx-\int_{\Omega}{\Dfdelta {l+n}}(\nabla_\epsilon u)X_{2n+1}(\eta^2w^\beta X^\epsilon_lu)dx
\end{equation*}
Expanding the first term we obtain
\begin{equation*}
\aligned
	0 = &-\sum_{i,j=1}^{2n+1} \int_\Omega{\Dfdeltaij} (\nabla_\epsilon u)X^\epsilon_lX^\epsilon_juX^\epsilon_lX^\epsilon_i u \eta^2w^\beta dx\\&+\sum_{j=1}^{2n+1}\int_\Omega{\Dfdeltajnl} (\nabla_\epsilon u)X^\epsilon_lX^\epsilon_juX_{2n+1}u \eta^2w^\beta dx\\&-\sum_{i,j=1}^{2n+1}\int_\Omega{\Dfdeltaij} (\nabla_\epsilon u)X^\epsilon_lX^\epsilon_juX_i^\epsilon(\eta^2w^\beta)  X^\epsilon_ludx\\
	&-\int_{\Omega}{\Dfdelta {l+n}}(\nabla_\epsilon u)X_{2n+1}(\eta^2w^\beta X^\epsilon_lu)dx.
\endaligned
\end{equation*}
Hence, bringing to the left hand side the first and the third term, we obtain
\begin{equation}\label{aiutiamoilettore}
\aligned
    &\int_\Omega\eta^2w^\beta\sum_{i,j=1}^{2n+1}{\Dfdeltaij} (\nabla_\epsilon u)X^\epsilon_lX^\epsilon_juX^\epsilon_lX_i^\epsilon u  dx\\&+\int_\Omega X^\epsilon_lu\eta^2\sum_{i,j=1}^{2n+1}{\Dfdeltaij} (\nabla_\epsilon u)X^\epsilon_lX^\epsilon_juX_i^\epsilon(w^\beta) dx\\&=-2\int_\Omega \eta X^\epsilon_luw^\beta\sum_{i,j=1}^{2n+1}{\Dfdeltaij} (\nabla_\epsilon u)X^\epsilon_lX^\epsilon_juX_i^\epsilon\eta dx\\&+\int_\Omega\eta^2w^\beta\sum_{j=1}^{2n+1}{\Dfdeltajnl} (\nabla_\epsilon u)X^\epsilon_lX^\epsilon_juX_{2n+1}udx\\&-\int_{\Omega}{\Dfdelta {l+n}}(\nabla_\epsilon u)X_{2n+1}(\eta^2w^\beta X^\epsilon_lu)dx=I_1^l+I_2^l+I_3^l.
\endaligned
\end{equation}
If $l\in\{n+1,\ldots,2n\}$ we argue in the same way starting from equation \eqref{weakform} and obtain
\begin{equation}\label{facciamocicapire}
\aligned
    &\int_\Omega\eta^2w^\beta\sum_{i,j=1}^{2n+1}{\Dfdeltaij} (\nabla_\epsilon u)X^\epsilon_lX^\epsilon_juX^\epsilon_lX_i^\epsilon udx\\&+\int_\Omega X^\epsilon_lu\eta^2\sum_{i,j=1}^{2n+1}{\Dfdeltaij} (\nabla_\epsilon u)X^\epsilon_lX^\epsilon_juX_i^\epsilon(w^\beta)dx\\&=-2\int_\Omega \eta X^\epsilon_luw^\beta\sum_{i,j=1}^{2n+1}{\Dfdeltaij} (\nabla_\epsilon u)X^\epsilon_lX^\epsilon_juX_i^\epsilon\eta dx\\&-\int_\Omega\eta^2w^\beta\sum_{j=1}^{2n+1}{\Dfdeltajln}(\nabla_\epsilon u)X^\epsilon_lX^\epsilon_juX_{2n+1}udx\\&+ \int_{\Omega}{\Dfdelta {l-n}}(\nabla_\epsilon u)X_{2n+1}(\eta^2w^\beta X^\epsilon_lu)dx=I_1^l+I_2^l+I_3^l.
\endaligned
\end{equation}
In the end, if $l=2n+1$, we start again from equation \eqref{weakform} and using the fact that $[X_{2n+1}^\epsilon,X^\epsilon_i]=0$, for any $i=1,\ldots,2n+1$, implies that $I_2^{2n+1}=0$ and $I_3^{2n+1}=0$, and we deduce
\begin{equation}\label{numerare}
\aligned
    &\int_\Omega\eta^2w^\beta\sum_{i,j=1}^{2n+1}{\Dfdeltaij} (\nabla_\epsilon u)X^\epsilon_{l}X^\epsilon_juX^\epsilon_{l}X_i^\epsilon udx\\&+\int_\Omega X^\epsilon_{l}u\eta^2\sum_{i,j=1}^{2n+1}{\Dfdeltaij} (\nabla_\epsilon u)X^\epsilon_{l}X^\epsilon_juX_i^\epsilon(w^\beta)dx\\&=-2\int_\Omega \eta X^\epsilon_{l}uw^\beta\sum_{i,j=1}^{2n+1}{\Dfdeltaij} (\nabla_\epsilon u)X^\epsilon_{l}X^\epsilon_juX_i^\epsilon\eta dx=:I_1^{2n+1}.
\endaligned
\end{equation}

We have found the three equalities \eqref{aiutiamoilettore}, \eqref{facciamocicapire} and \eqref{numerare} for different values of $l=1,\ldots,2n+1$. Since the left hand side is the same for all values of $l,$
we will handle together  the second integral in the left hand side, for any $l=1,\ldots,2n+1$:
\begin{equation*}
\aligned
	&\int_\Omega X^\epsilon_lu\eta^2\sum_{i,j=1}^{2n+1}{\Dfdeltaij} (\nabla_\epsilon u)X^\epsilon_lX^\epsilon_juX_i^\epsilon(w^\beta)dx\\&=\beta\int_\Omega \eta^2w^{\beta-2}\sum_{i,j=1}^{2n+1}{\Dfdeltaij} (\nabla_\epsilon u)X^\epsilon_lX^\epsilon_juX^\epsilon_lu\sum_{k=1}^{2n+1}X_i^\epsilon X^\epsilon_kuX^\epsilon_kudx\\&=\beta\int_\Omega \eta^2w^{\beta-2}\sum_{i,j=1}^{2n+1}{\Dfdeltaij} (\nabla_\epsilon u)X^\epsilon_lX^\epsilon_juX^\epsilon_lu\sum_{k=1}^{2n+1}X^\epsilon_kX_i^\epsilon uX^\epsilon_kudx\\&-\beta\sigma(i)\int_\Omega \eta^2w^{\beta-2}X^\epsilon_lu X_{2n+1}u\sum_{i,j=1}^{2n} {\Dfdeltaij} (\nabla_\epsilon u)X^\epsilon_lX^\epsilon_ju X^\epsilon_{i-\sigma(i)n}udx,
\endaligned
\end{equation*}
where we called $\sigma(k)=-1$ if $k\leq n$, $\sigma(k)=1$ if $k> n$, so that $[X^\epsilon_i,X^\epsilon_k]=\sigma(k)\delta_{i,k-\sigma(k) n} X_{2n+1}$ for any $i,k\in\left\{1,\ldots,2n\right\}$.
   
Denoting by $I^l_4$ the last integral, summing up over $l$ the three equalities \eqref{aiutiamoilettore}, \eqref{facciamocicapire} and \eqref{numerare} and summing everything together we get
\begin{align}\label{semprenumerare}
	&\int_\Omega\eta^2w^\beta\sum_{i,j,l=1}^{2n+1}{\Dfdeltaij} (\nabla_\epsilon u)X^\epsilon_lX^\epsilon_juX^\epsilon_lX_i^\epsilon udx\\&+ \beta \int_\Omega \eta^2w^{\beta-2}\sum_{i,j=1}^{2n+1}{\Dfdeltaij} (\nabla_\epsilon u)g_\epsilon\left(\nabla_\epsilon X^\epsilon_j,\nabla_\epsilon u\right) g_\epsilon\left(\nabla_\epsilon X_i^\epsilon u,\nabla_\epsilon u\right) dx\nonumber \\&=
	\sum_{l=1}^{2n+1}\left(I_1^l+I_2^l+I_3^l-I_4^l\right).\nonumber
\end{align}
The second integral in the left hand side of \eqref{semprenumerare} is always positive, hence, from the structure condition \eqref{structurecondition01} and using the notations \eqref{Phi_def} and  \eqref{26ottobre},  \eqref{semprenumerare} reduces to
\begin{equation}\label{numeriamoancora}
	\int_\Omega\eta^2 w^\beta \|\Phi\|_\lambda^2 dx\leq\sum_{l=1}^{2n+1}\left(I_1^l+I_2^l+I_3^l+|I_4^l|\right),
\end{equation}

where $I_1^l$ have been defined in equations 
\eqref{aiutiamoilettore}, \eqref{facciamocicapire} and \eqref{numerare} respectively. 
 
Moreover, from \eqref{structurecondition01} and Young's inequality, again using the notations \eqref{Phi_def} and \eqref{26ottobre}, it follows that
\begin{align*}
	|I_1^l|&\leq 2\int_\Omega \eta w^{\beta+1}\sum_{i,j=1}^{2n+1}{\Dfdeltaij} (\nabla_\epsilon u)X^\epsilon_lX^\epsilon_juX_i^\epsilon\eta dx
	\\&\leq c\tau\int_\Omega \eta^2w^\beta\sum_{i=1}^{n}\lambda_{i,\delta}(\nabla_\epsilon u)\left(|\nabla_\epsilon X_i^\epsilon u|_\epsilon^2+|\nabla_\epsilon X^\epsilon_{i+n}u|_\epsilon^2+|\nabla_\epsilon X^\epsilon_{2n+1}u|_\epsilon^2\right)dx\\&+\frac{c}{\tau}\int_\Omega w^{\beta+2}\sum_{i=1}^{n}\lambda_{i,\delta}(\nabla_\epsilon u)|\nabla_\epsilon\eta|_\epsilon^2 dx\\&=c\tau\int_\Omega\eta^2w^\beta\|\Phi\|_\lambda^2 dx+\frac{c}{\tau}\int_\Omega|\nabla_\epsilon \eta|_\epsilon^2w^{\beta+2}\Lambda dx,
\end{align*}
for all $l\in\{1,\ldots,2n+1\}$, where $c=c(q,n,L)>0$. For $\tau$ small enough the first term can be reabsorbed in the left hand side of equation \eqref{numeriamoancora}.

Let us estimate $I_2^l$, for $l\in\{1,\ldots,2n\}$: again by \eqref{structurecondition01} and Young's inequality, using the notations \eqref{Phi_def} and \eqref{26ottobre}, we have
\begin{align*}
	|I_2^l|&\leq c\tau \int_\Omega\eta^2w^\beta\sum_{i=1}^{n}\lambda_{i,\delta}(\nabla_\epsilon u)\left(|\nabla_\epsilon X_i^\epsilon u|_\epsilon^2+|\nabla_\epsilon X^\epsilon_{i+n}u|_\epsilon^2+|\nabla_\epsilon X^\epsilon_{2n+1}u|_\epsilon^2\right)dx\\&+\frac{c}{\tau}\int_\Omega \eta^2w^\beta\sum_{i=1}^{n}\lambda_{i,\delta}(\nabla_\epsilon u)|X_{2n+1}u|^2dx\\
	&=c\tau\int_\Omega\eta^2w^\beta\|\Phi\|_\lambda^2 dx+\frac{c}{\tau}\int_\Omega\eta^2w^\beta\Lambda|X_{2n+1}u|^2dx.
\end{align*}
As before, for $\tau$ small enough the first term can be reabsorbed in the left hand side of equation \eqref{numeriamoancora}.
In order to estimate $I_3^l$, for all $l\in\{1,\ldots,2n\}$, we explicitly compute the derivative in the integrand 
\begin{equation*}
	X_{2n+1}(\eta^2w^\beta X^\epsilon_lu)=2\eta X_{2n+1}\eta w^\beta X^\epsilon_lu+\beta\eta^2X^\epsilon_luw^{\beta-2}\sum_{k=1}^{2n+1}X^\epsilon_kuX^\epsilon_kX_{2n+1}u+\eta^2w^\beta X^\epsilon_lX_{2n+1}u.
\end{equation*}
If $l\in\{1,\ldots,n\}$ 
\begin{align*}
	I_3^l&=-2\int_{\Omega}\eta X_{2n+1}\eta {\Dfdelta {l+n}}(\nabla_\epsilon u)w^\beta X^\epsilon_ludx-\int_{\Omega}\eta^2w^\beta{\Dfdelta {l+n}}(\nabla_\epsilon u)X^\epsilon_lX_{2n+1}udx \\
	&-\beta\int_{\Omega}\eta^2X^\epsilon_luw^{\beta-2}{\Dfdelta {l+n}}(\nabla_\epsilon u)\sum_{k=1}^{2n+1}X^\epsilon_kuX^\epsilon_kX_{2n+1}udx=I_{3,1}^l+I_{3,2}^l+I_{3,3}^l,
\end{align*}
otherwise if $l\in\{n+1,\ldots,2n\}$
\begin{align*}
	I_3^l&=2\int_{\Omega}\eta X_{2n+1}\eta {\Dfdelta {l-n}}(\nabla_\epsilon u)w^\beta X^\epsilon_ludx+\int_{\Omega}\eta^2w^\beta{\Dfdelta {l-n}}(\nabla_\epsilon u)X^\epsilon_lX_{2n+1}udx \\
	&+\beta\int_{\Omega}\eta^2X^\epsilon_luw^{\beta-2}{\Dfdelta {l-n}}(\nabla_\epsilon u)\sum_{k=1}^{2n+1}X^\epsilon_kuX^\epsilon_kX_{2n+1}udx=I_{3,1}^l+I_{3,2}^l+I_{3,3}^l,
\end{align*}
From \eqref{structurecondition02} and \eqref{26ottobre} it follows that 
\begin{equation*}
	|I_{3,1}^l|\leq c\int_\Omega \eta|X_{2n+1}\eta|w^\beta\lambda_{l,\delta}(\nabla_\epsilon u)|X_l^\epsilon u| |X_{l+n}^\epsilon u|dx\leq c\int_\Omega \eta|X_{2n+1}\eta|w^{\beta+2}\Lambda dx,
\end{equation*}
for any $l\in\{1,\ldots,n\}$; analogously if $l\in\{n+1,\ldots,2n\}$
\begin{equation*}
	|I_{3,1}^l|\leq c\int_\Omega \eta|X_{2n+1}\eta|w^{\beta+2}\Lambda dx.
\end{equation*}
Hence
\begin{equation*}
	|I_{3,1}^l|\leq c\int_\Omega \eta|X_{2n+1}\eta|w^{\beta+2}\Lambda dx,
\end{equation*}
for any $l\in\{1,\ldots,2n\}$, where $c=c(q,n,L)$.
		
In order to estimate $I_{3,2}^l$, we integrate by parts and we compute the derivative of the integrand function: hence for $l\in\{1,\ldots,n\}$
\begin{equation}\label{I32}
\aligned
	I_{3,2}^l&= \int_\Omega X^\epsilon_l(\eta^2w^\beta {\Dfdelta {l+n}}(\nabla_\epsilon u))X_{2n+1}udx=
	\int_\Omega \eta^2w^\beta\sum_{j=1}^{2n+1}{\Dfdeltajnl} (\nabla_\epsilon u)X^\epsilon_lX^\epsilon_juX_{2n+1}udx\\&+2\int_\Omega \eta X^\epsilon_l\eta {\Dfdelta {l+n}}(\nabla_\epsilon u)w^\beta X_{2n+1}udx
	+\beta\int_\Omega \eta^2{\Dfdelta {l+n}}(\nabla_\epsilon u)w^{\beta-2}\sum_{k=1}^{2n+1}X^\epsilon_lX^\epsilon_kuX^\epsilon_kuX_{2n+1}udx.
\endaligned
\end{equation}
	
The first integral in \eqref{I32} can estimated as $I_2^l$, hence we will estimate the last two of them: again by using \eqref{structurecondition02}, Young's inequality and \eqref{26ottobre}, the second term in \eqref{I32} can be bounded by
\begin{equation*}
\aligned
	&c\int_\Omega \eta|X^\epsilon_l\eta| \lambda_{l,\delta}(\nabla_\epsilon u)w^\beta|X_{l+n}^\epsilon u| |X_{2n+1}u|dx\\&\leq c\int_\Omega \eta|\nabla_\epsilon\eta|_\epsilon \sum_{i=1}^{n}\lambda_{i,\delta}(\nabla_\epsilon u)w^{\beta+1}|X_{2n+1}u|dx\\&\leq c\int_\Omega|\nabla_\epsilon\eta|_\epsilon^2w^{\beta+2}\Lambda dx+c\int_\Omega\eta^2w^\beta\Lambda|X_{2n+1}u|^2dx,
\endaligned
\end{equation*}
for some constant $c=c\left(q,n,L\right)$. Instead, for the third one we use $[X_l^\epsilon,X_k^\epsilon]=X_{2n+1}$ if $k=l+n$, then 
\begin{equation}\label{4gennaio}
\aligned
	&\beta\int_\Omega \eta^2{\Dfdelta {l+n}}(\nabla_\epsilon u)w^{\beta-2}\sum_{k=1}^{2n+1}X^\epsilon_lX^\epsilon_kuX^\epsilon_kuX_{2n+1}udx\\&=\beta\int_\Omega \eta^2{\Dfdelta {l+n}}(\nabla_\epsilon u)w^{\beta-2}\sum_{k=1}^{2n+1}X^\epsilon_kX^\epsilon_luX^\epsilon_kuX_{2n+1}udx\\&+\beta\int_\Omega \eta^2{\Dfdelta {l+n}}(\nabla_\epsilon u)w^{\beta-2}X^\epsilon_{l+n}u|X_{2n+1}u|^2dx
\endaligned
\end{equation}
(by \eqref{structurecondition02}, Young's inequality and \eqref{26ottobre})
\begin{equation*}
\aligned
	&\leq c\beta\tau\int_\Omega \eta^2w^\beta \lambda_{l,\delta}(\nabla_\epsilon u)|\nabla_\epsilon X_l^\epsilon u|_\epsilon^2dx\\&+\frac{c\beta}{\tau}\int_\Omega \eta^2w^\beta\Lambda |X_{2n+1}u|^2dx+c\beta\int_\Omega\eta^2w^\beta\Lambda|X_{2n+1}u|^2dx\\& 
	\leq c\beta\tau\int_\Omega \eta^2w^\beta \lambda_{l,\delta}(\nabla_\epsilon u)\left(|\nabla_\epsilon X_l^\epsilon u|_\epsilon^2+|\nabla_\epsilon X_{l+n}^\epsilon u|_\epsilon^2+|\nabla_\epsilon X_{2n+1}^\epsilon u|_\epsilon^2\right)dx\\&+c\beta\left(1+\frac{1}{\tau}\right)\int_\Omega\eta^2w^\beta\Lambda|X_{2n+1}u|^2dx,
\endaligned
\end{equation*}
where again $c=c(q,n,L)$.
	
If $l\in\{n+1,\ldots,2n\}$ the same bounds hold for $I_{3,2}^l$. Therefore, using the notation \eqref{Phi_def},
\begin{equation*}
	|\sum_{l=1}^{2n}I_{3,2}^l|\leq c\beta\tau\int_\Omega \eta^2w^\beta\|\Phi\|_\lambda^2 dx+c\beta\left(1+\frac{1}{\tau}\right)\int_\Omega\eta^2w^\beta|X_{2n+1}u|^2\Lambda dx.
\end{equation*}
The integral $I_{3,3}^l$ can be estimated as \eqref{4gennaio}.
	
In the end, for $I^l_4$ again \eqref{structurecondition02}, Young's inequality and \eqref{26ottobre} imply
\begin{equation*}
\aligned
	|I_4^l|&\leq\beta\int_\Omega \eta^2w^{\beta-1} |X_{2n+1}u|\\&\left(\sum_{i=1}^{n}\lambda_{i,\delta}(\nabla_\epsilon u)\left( (X^\epsilon_lX_i^\epsilon u)^2+(X^\epsilon_lX_{i+n}^\epsilon u)^2+(X^\epsilon_lX_{2n+1}^\epsilon u)^2\right)\right)^{\frac{1}{2}}\left(\sum_{i=1}^{n}\lambda_{i,\delta}(\nabla_\epsilon u)|\nabla_\epsilon u|_\epsilon^2\right)^{\frac{1}{2}}dx\\&=c\beta\tau\int_\Omega\eta^2w^\beta\|\Phi\|_\lambda^2 dx+\frac{c\beta}{\tau}\int_\Omega\eta^2w^\beta|X_{2n+1}u|^2 \Lambda dx.
\endaligned
\end{equation*}
Hence the thesis follows choosing $\tau=\frac{1}{2c(\beta+1)}$.
\end{proof}

Let us note explicitly that Lemma \ref{caccioppoliXuriemannian} cannot be considered a real Caccioppoli type inequality, due to the presence of the term
\begin{equation}\label{badterm}
	\int_\Omega \eta^2w^{\beta}(\nabla_\epsilon u)|X_{2n+1}u|^2\Lambda dx.
\end{equation}
The main goal is now to obtain a better inequality, without this term. The first step is to prove the following Caccioppoli-type inequality for $X_{2n+1}u$.

\begin{lemma}[Caccioppoli for $X_{2n+1}u$]\label{eq:CaccioppoliTuapprox}
There exists a constant $c=c(q,n,L)>0$, independent of $\epsilon$ and $\delta$, such that, for every weak solution $u\in W_\epsilon^{1,q}(\Omega)$ to \eqref{eqneriemannianapprox}, for every $\beta\geq0$ and for every non-negative $\eta\in C_c^\infty(\Omega)$ one has
\begin{equation*}
    \int_\Omega \eta^2\ |X_{2n+1}u|^\beta\ \|\Psi(X_{2n+1}u)\|^2_\lambda\  dx\leq\frac{c}{(\beta+1)^2} \int_\Omega |\nabla_\epsilon\eta|_\epsilon^2\ |X_{2n+1}u|^{\beta+2}\ \Lambda dx,
\end{equation*}
where $\Psi$ is defined in \eqref{26ottobre1} and $\Lambda$ in \eqref{26ottobre}.
\end{lemma}

\begin{proof}
Let us consider $\varphi=X_{2n+1}(\eta^2|X_{2n+1}u|^\beta X_{2n+1}u)$ as test function in the weak formulation \eqref{weakformulation} of the equation \eqref{eqneriemannianapprox}
\begin{equation*}
	\sum_{i=1}^{2n+1}\int_\Omega {\Dfdelta i}(\nabla_\epsilon u)X_i^\epsilon(X_{2n+1}(\eta^2|X_{2n+1}u|^\beta X_{2n+1}u))dx=0.
\end{equation*}
Since the vector field $X_{2n+1}$ commutes with $X_i^\epsilon$, for any $i=1,\ldots,2n+1$, an integration by parts gives
\begin{equation*}
\aligned
    &\int_\Omega\eta^2|X_{2n+1}u|^\beta\sum_{i,j=1}^{2n+1}{\Dfdeltaij}(\nabla_\epsilon u)X_j^\epsilon X_{2n+1}u X_i^\epsilon X_{2n+1}u dx\\&=-\int_{\Omega}X_{2n+1}u\sum_{i,j=1}^{2n+1}{\Dfdeltaij}(\nabla_\epsilon u)X_j^\epsilon X_{2n+1}u X_i^\epsilon (\eta^2|X_{2n+1}u|^\beta)dx.
\endaligned
\end{equation*}
By derivating on the right hand side, one gets
\begin{equation*}
\aligned
    &\int_\Omega\eta^2|X_{2n+1}u|^\beta\sum_{i,j=1}^{2n+1}{\Dfdeltaij}(\nabla_\epsilon u)X_j^\epsilon X_{2n+1}u X_i^\epsilon X_{2n+1}u dx\\&\leq \frac{2}{(\beta+1)}\int_{\Omega}\eta|X_{2n+1}u|^{\beta+1}\bigg{\vert}\sum_{i,j=1}^{2n+1}{\Dfdeltaij}(\nabla_\epsilon u)X_j^\epsilon X_{2n+1}u X_i^\epsilon \eta\bigg{\vert} dx.
\endaligned
\end{equation*}
From the structure conditions \eqref{structurecondition01}, and using the notation \eqref{26ottobre1}, it follows that the left hand side
\begin{equation*}   (LHS)\geq\int_\Omega \eta^2|X_{2n+1}u|^\beta\|\Psi\|_\lambda^2 dx.
\end{equation*}
As for the right hand side, by \eqref{structurecondition01} and Young's inequality, using the notation \eqref{26ottobre}, it follows
\begin{equation*}
	(RHS)\leq\frac{\tau}{2} \int_\Omega\eta^2|X_{2n+1}u|^\beta\|\Psi\|_\lambda^2 dx+\frac{2L^2}{\tau(\beta+1)^2} \int_\Omega|\nabla_\epsilon\eta|_\epsilon^2|X_{2n+1}u|^{\beta+2}\Lambda dx.
\end{equation*}
Choosing $\tau=1$, we obtain
\begin{equation*}
    \int_\Omega \eta^2|X_{2n+1}u|^\beta\|\Psi\|_\lambda^2 dx\leq \frac{c}{(\beta+1)^2}\int_\Omega|\nabla_\epsilon\eta|_\epsilon^2 |X_{2n+1}u|^{\beta+2}\Lambda dx,
\end{equation*}
with $c=4L^2$ and this finishes the proof.
\end{proof}

The mixed Caccioppoli-type inequality in this case reads as follows.

\begin{lemma}\label{mainlemma}
There exists a constant $c=c(q,n,L)>0$, independent of $\epsilon$ and $\delta$, such that, for every weak solution
$u\in W_\epsilon^{1,q}(\Omega)$ to \eqref{eqneriemannianapprox}, for every $\alpha,\beta\geq2$ and for every non-negative $\eta\in C_c^\infty(\Omega)$ one has
\begin{equation*}
    \int_\Omega \eta^{\alpha+2}\ |X_{2n+1}u|^{\beta}\ \|\Phi(\nabla_\epsilon u)\|_\lambda^2\ dx\leq c(\alpha+2)^2\|\nabla_\epsilon\eta\|_\infty^2\int_\Omega \eta^\alpha\ w(\nabla_\epsilon u)^2\ |X_{2n+1}u|^{\beta-2}\ \|\Phi(\nabla_\epsilon u)\|_\lambda^2\ dx,
\end{equation*}
where $\|\Phi\|_\lambda^2$ is defined in \eqref{Phi_def} and $w$ is defined in \eqref{26ottobre}.
\end{lemma}

\begin{proof}
Let $\eta\in C_c^\infty(\Omega)$ be a non-negative cut-off function, $l\in\{1,\ldots,2n+1\}$ and consider $\varphi=X_l^\epsilon(\eta^{\alpha+2}|X_{2n+1}u|^\beta X_l^\epsilon u )$ as test function in the weak formulation \eqref{weakformulation} of the equation \eqref{eqneriemannianapprox}. Let us compute the derivative
\begin{equation*}
\aligned
    X_i^\epsilon(\eta^{\alpha+2}|X_{2n+1}u|^\beta X_l^\epsilon u)&=(\alpha+2)\eta^{\alpha+1} X_i^\epsilon\eta|X_{2n+1}u|^\beta X_l^\epsilon u\\&+\beta\eta^{\alpha+2}|X_{2n+1}u|^{\beta-2}X_{2n+1}uX_l^\epsilon uX_i^\epsilon X_{2n+1}u+\eta^{\alpha+2}|X_{2n+1}u|^\beta X_i^\epsilon X_l^\epsilon u.
\endaligned
\end{equation*}
Let us suppose that $l\in\{1,\ldots,n\}$, we obtain 
\begin{equation}\label{30dicembre}
\aligned
	&\sum_{i=1}^{2n+1}\int_\Omega X_l^\epsilon({\Dfdelta i}(\nabla_\epsilon u))X^\epsilon_lX^\epsilon_iu 
	|X_{2n+1}u|^{\beta}\eta^{\alpha+2} dx\\&=\int_\Omega X_l^\epsilon({\Dfdelta {l+n}}(\nabla_\epsilon u))X_{2n+1}u |X_{2n+1}u|^{\beta}\eta^{\alpha+2}dx\\
	&-(\alpha+2)\sum_{i=1}^{2n+1}\int_\Omega X^\epsilon_l({\Dfdelta i}(\nabla_\epsilon u))X^\epsilon_lu |X_{2n+1}u|^{\beta}\eta^{\alpha+1}X^\epsilon_i\eta dx\\
	&+\int_\Omega X_{2n+1}({\Dfdelta {l+n}}(\nabla_\epsilon u))|X_{2n+1}u|^{\beta}X^\epsilon_lu \eta^{\alpha+2}dx\\
	&-\beta\sum_{i=1}^{2n+1}\int_\Omega X^\epsilon_l({\Dfdelta i}(\nabla_\epsilon u)X^\epsilon_iX_{2n+1}u X^\epsilon_lu|X_{2n+1}u|^{\beta-2}X_{2n+1}u \eta^{\alpha+2}dx=:I_1^l+I_2^l+I_3^l+I_4^l.
\endaligned
\end{equation}
	
From the structure conditions \eqref{structurecondition01} it follows that 
\begin{equation*}
	(LHS)\geq \int_\Omega \eta^{\alpha+2}|X_{2n+1}u|^{\beta}\sum_{i=1}^n\lambda_{i,\delta}(\nabla_\epsilon u)((X^\epsilon_lX^\epsilon_iu)^2+(X^\epsilon_lX^\epsilon_{i+n}u)^2+(X_l^\epsilon X_{2n+1}^\epsilon u)^2)dx.
\end{equation*}
We will prove that
\begin{equation}\label{29dicembre3}
	|I_k^l|\leq c\tau \int_\Omega \eta^{\alpha+2}|X_{2n+1}u|^{\beta}\|\Phi\|_\lambda^2 dx
	+\frac{c(\alpha+2)^2\|\nabla_\epsilon\eta\|_\infty^2}{\tau}\int_\Omega \eta^\alpha w^2|X_{2n+1}u|^{\beta-2}\|\Phi\|_\lambda^2 dx,
\end{equation}
for some $c=c(q,n,L)>0$, for any $k=1,2,3,4$.
	
Let us start by estimating $I_4^l$: via the Cauchy-Schwarz inequality and the structure condition \eqref{structurecondition01}, and using the notation \eqref{Phi_def} and \eqref{26ottobre}, we have
\begin{align*}
	|I_4^l| &\leq \frac{c\tau(\beta+1)^2}{\|\nabla_\epsilon\eta\|_\infty^2(\alpha+2)^2}\int_\Omega \eta^{\alpha+4} |X_{2n+1}u|^{\beta}\|\Psi\|_\lambda^2 dx+\frac{c(\alpha+2)^2\|\nabla_\epsilon\eta\|_\infty^2}{\tau}\\&\int_\Omega\eta^{\alpha}|X_{2n+1}u|^{\beta-2}w^2\sum_{i=1}^{n}\lambda_{i,\delta}(\nabla_\epsilon u)\left((X^\epsilon_lX^\epsilon_iu)^2+(X^\epsilon_lX^\epsilon_{i+n}u)^2+(X_l^\epsilon uX_{2n+1}^\epsilon u)^2\right) dx,
\end{align*}
where $c=c(q,n,L)>0$.

By  Lemma \ref{eq:CaccioppoliTuapprox} the first term can be bounded by
\begin{equation}\label{estimate1}
	\tau c \int_\Omega \eta^{\alpha+2}|X_{2n+1}u|^{\beta}|X_{2n+1}u|^2\Lambda dx.
\end{equation}
Let us observe that, by definition
\begin{equation}\label{29dicembre}
\aligned
    \Lambda |X_{2n+1}|^2&=\sum_{i=1}^{n}\lambda_{i,\delta}(\nabla_\epsilon u)|X_{2n+1}u|^{2}=\sum_{i=1}^{n}\lambda_{i,\delta}(\nabla_\epsilon u)|X^\epsilon_iX^\epsilon_{i+n}u-X^\epsilon_{i+n}X^\epsilon_iu|^{2}\\ &\leq \sum_{i=1}^{n}\lambda_{i,\delta}(\nabla_\epsilon u)\left(|\nabla_\epsilon X_i^\epsilon u|_\epsilon^2+|\nabla_\epsilon X_{i+n}^\epsilon u|_\epsilon^2\right)\\
	&\leq \sum_{i=1}^{n}\lambda_{i,\delta}(\nabla_\epsilon u)\left(|\nabla_\epsilon X_i^\epsilon u|_\epsilon^2+|\nabla_\epsilon X_{i+n}^\epsilon u|_\epsilon^2+|\nabla_\epsilon X_{2n+1}^\epsilon u|_\epsilon^2\right)=\|\Phi\|_\lambda^2,
\endaligned
\end{equation}
where in the last inequality we simply added a positive term.
	
Hence, \eqref{estimate1} can be bounded by
\begin{equation}\label{29dicembre2}
	c\tau \int_\Omega \eta^{\alpha+2}|X_{2n+1}u|^{\beta}\|\Phi\|_\lambda^2 dx,
\end{equation}
which will be reabsorbed in the left hand side for $\tau$ small enough.
		
Let us estimate $I_3^l$:
\begin{equation*}
\aligned
	I_3^l&=\int_\Omega \sum_{j=1}^{2n}{\Dfdeltajnl} (\nabla_\epsilon u)X^\epsilon_jX_{2n+1}u|X_{2n+1}u|^{\beta}X^\epsilon_lu\eta^{\alpha+2}\ dx\\
	&\leq \frac{c\tau(\beta+1)^2}{\|\nabla_\epsilon\eta\|_\infty^2(\alpha+2)^2}\int_\Omega \eta^{\alpha+4}|X_{2n+1}u|^{\beta}\|\Psi\|_\lambda^2dx\\&+\frac{c\|\nabla_\epsilon\eta\|_\infty^2(\alpha+2)^2}{\tau(\beta+1)^2}\int_\Omega\eta^{\alpha}|X_{2n+1}u|^{\beta}w^2\lambda_{l,\delta}(\nabla_\epsilon u) dx,
\endaligned
\end{equation*}
where $c=c(q,n,L)>0$.
The first term can be handled by using Lemma \ref{eq:CaccioppoliTuapprox} and \eqref{29dicembre}, hence it can be bounded by \eqref{29dicembre2}; while for the second we use \eqref{29dicembre}
\begin{equation*}
\aligned
    &\frac{c\|\nabla_\epsilon\eta\|_\infty^2(\alpha+2)^2}{\tau(\beta+1)^2}\int_\Omega\eta^{\alpha}|X_{2n+1}u|^{\beta-2}w^2\lambda_{l,\delta}(\nabla_\epsilon u)|X_{2n+1}u|^2 dx\\
	&\leq\frac{c\|\nabla_\epsilon\eta\|_\infty^2(\alpha+2)^2}{\tau(\beta+1)^2}\int_\Omega\eta^{\alpha}|X_{2n+1}u|^{\beta-2}w^2\Lambda|X_{2n+1}u|^2dx\\&\leq
	\frac{c(\alpha+2)^2\|\nabla_\epsilon\eta\|_\infty^2}{\tau}\int_\Omega\eta^{\alpha}|X_{2n+1}u|^{\beta-2}w^2\|\Phi\|_\lambda^2dx.
\endaligned
\end{equation*}
Let us estimate
$$I_2^l=-(\alpha+2)\int_\Omega |X_{2n+1}u|^{\beta}\eta^{\alpha+1}\sum_{i,j=1}^{2n+1}{\Dfdeltaij}(\nabla_\epsilon u)X^\epsilon_lX^\epsilon_juX^\epsilon_lu X^\epsilon_i\eta dx.$$
Again by the structure condition \eqref{structurecondition01}, Cauchy-Schwartz's and Young's inequalities
\begin{equation*}
\aligned
	|I_2^l|&\leq
	\frac{c\tau}{\|\nabla_\epsilon\eta\|_\infty^2}\int_\Omega \eta^{\alpha+2}|X_{2n+1}u|^{\beta+2}\sum_{i=1}^{n}\lambda_{i,\delta}(\nabla_\epsilon u)((X^\epsilon_i\eta)^2+(X^\epsilon_{i+n}\eta)^2+(X^\epsilon_{2n+1}\eta)^2)dx\\&+\frac{c(\alpha+2)^2\|\nabla_\epsilon\eta\|_\infty^2}{\tau}\int_\Omega \eta^\alpha w^2|X_{2n+1}u|^{\beta-2}\|\Phi\|_\lambda^2dx.
\endaligned
\end{equation*}
By using \eqref{29dicembre} the first term can be bounded by
\begin{equation*}
	c\tau \int_\Omega \eta^{\alpha+2}|X_{2n+1}u|^{\beta}\|\Phi\|_\lambda^2dx.
\end{equation*}
In the end
\begin{equation*}
\aligned
	I_1^l&=-(\beta+1)\int_\Omega \eta^{\alpha+2}|X_{2n+1}u|^{\beta}{\Dfdelta {l+n}}(\nabla_\epsilon u) X^\epsilon_lX_{2n+1}u dx\\
	&-(\alpha+2)\int_\Omega \eta^{\alpha+1}|X_{2n+1}u|^{\beta}X_{2n+1}u{\Dfdelta {l+n}}(\nabla_\epsilon u)X^\epsilon_l\eta dx=I_{1,1}^l+I_{1,2}^l.
\endaligned
\end{equation*}

As for $I_{1,1}^l$, \eqref{structurecondition02}, the Young inequality, \eqref{29dicembre} and Lemma \ref{eq:CaccioppoliTuapprox} imply that
\begin{equation*}
\aligned
	|I_{1,1}^l|&\leq (\beta+1)\int_\Omega\eta^{\alpha+2}|X_{2n+1}u|^\beta |X^\epsilon_lX_{2n+1}u|\lambda_{l,\delta}(\nabla_\epsilon u)|X_l^\epsilon u| dx\\
	&\leq\frac{c\tau(\beta+1)^2}{\|\nabla_\epsilon\eta\|_\infty^2(\alpha+2)^2}\int_\Omega\eta^{\alpha+4}|X_{2n+1}u|^\beta\lambda_{l,\delta}(\nabla_\epsilon u)(X^\epsilon_lX_{2n+1}u)^2dx\\       &+\frac{c(\alpha+2)^2\|\nabla_\epsilon\eta\|_\infty^2}{\tau}\int_\Omega \eta^\alpha |X_{2n+1}u|^\beta w^2\lambda_{l,\delta}(\nabla_\epsilon u)dx\\
	&\leq\frac{c\tau(\beta+1)^2}{\|\nabla_\epsilon\eta\|_\infty^2(\alpha+2)^2}\int_\Omega\eta^{\alpha+4}|X_{2n+1}u|^\beta\|\Psi\|_\lambda^2dx\\
	&+\frac{c(\alpha+2)^2\|\nabla_\epsilon\eta\|_\infty^2}{\tau}\int_\Omega \eta^\alpha |X_{2n+1}u|^{\beta-2} w^2\Lambda|X_{2n+1}u|^2dx\\
	&\leq c\tau \int_\Omega \eta^{\alpha+2}|X_{2n+1}u|^{\beta}\|\Phi\|_\lambda^2dx
	+\frac{c(\alpha+2)^2\|\nabla_\epsilon\eta\|_\infty^2}{\tau}\int_\Omega \eta^\alpha w^2|X_{2n+1}u|^{\beta-2}\|\Phi\|_\lambda^2dx,
\endaligned
\end{equation*}
while for $I_{1,2}^l$, \eqref{structurecondition02}, the Young's inequality and \eqref{29dicembre} imply that
\begin{equation*}
\aligned
	|I_{1,2}^l|&\leq (\alpha+2)\int_\Omega \eta^{\alpha+1}\lambda_{l,\delta}(\nabla_\epsilon u)|\nabla_\epsilon u|_\epsilon|\nabla_\epsilon\eta|_\epsilon |X_{2n+1}u|^{\beta+1}dx\\
	&\leq\frac{c\tau}{\|\nabla_\epsilon\eta\|_\infty^2}\int_\Omega \eta^{\alpha+2}\lambda_{l,\delta}(\nabla_\epsilon u)|\nabla_\epsilon\eta|_\epsilon^2|X_{2n+1}u|^{\beta+2}dx\\&+
	\frac{c(\alpha+2)^2\|\nabla_\epsilon\eta\|_\infty^2}{\tau}\int_\Omega \eta^\alpha |X_{2n+1}u|^\beta w^2\lambda_{l,\delta}(\nabla_\epsilon u)dx\\&\leq
	c\tau\int_\Omega \eta^{\alpha+2}|X_{2n+1}u|^\beta\Lambda|X_{2n+1}u|^{2}dx\\&+
	\frac{c(\alpha+2)^2\|\nabla_\epsilon\eta\|_\infty^2}{\tau}\int_\Omega \eta^\alpha |X_{2n+1}u|^{\beta-2}w^2\Lambda|X_{2n+1}u|^2dx\\&\leq c\tau \int_\Omega \eta^{\alpha+2}|X_{2n+1}u|^{\beta}\|\Phi\|_\lambda^2dx+\frac{c(\alpha+2)^2\|\nabla_\epsilon\eta\|_\infty^2}{\tau}\int_\Omega \eta^\alpha w^2|X_{2n+1}u|^{\beta-2}\|\Phi\|_\lambda^2dx.
\endaligned
\end{equation*}
The bound \eqref{29dicembre3} holds true even if $l\in\{n+1,\ldots,2n\}$. 
	
If $l=2n+1$, the test function reads as $\varphi=X^\epsilon_{2n+1}(\eta^{\alpha+2}|X_{2n+1}u|^\beta X_{2n+1}^\epsilon u)$: since $[X_{2n+1}^\epsilon,X_i^\epsilon]=0, \forall i=1,\ldots,2n+1$, in \eqref{30dicembre} the terms $I_1^{2n+1}$ and $I_3^{2n+1}$ will disappear. Hence
\begin{align*}
	&\sum_{i=1}^{2n+1}\int_\Omega X_{2n+1}^\epsilon({\Dfdelta i}(\nabla_\epsilon u))X^\epsilon_{2n+1}X^\epsilon_iu 
	|X_{2n+1}u|^{\beta}\eta^{\alpha+2} dx\\&=-(\alpha+2)\sum_{i=1}^{2n+1}\int_\Omega X^\epsilon_{2n+1}({\Dfdelta i}(\nabla_\epsilon u))X^\epsilon_{2n+1}u |X_{2n+1}u|^{\beta}\eta^{\alpha+1}X^\epsilon_i\eta dx\\
	&-\beta\sum_{i=1}^{2n+1}\int_\Omega X^\epsilon_{2n+1}({\Dfdelta i}(\nabla_\epsilon u))X^\epsilon_iX_{2n+1}u X^\epsilon_{2n+1}u|X_{2n+1}u|^{\beta-2}X_{2n+1}u\eta^{\alpha+2}dx\\&=I_2^{2n+1}+I_4^{2n+1}.
\end{align*}
For $I_2^{2n+1}$ and $I_4^{2n+1}$ the bound \eqref{29dicembre3} holds. As for the left hand side, \eqref{structurecondition01} implies that
\begin{equation*}
	(LHS)\geq \int_\Omega \eta^{\alpha+2}|X_{2n+1}u|^{\beta}\sum_{i=1}^{n}\lambda_{i,\delta}(\nabla_\epsilon u)((X^\epsilon_{2n+1}X^\epsilon_iu)^2+(X^\epsilon_{2n+1}X^\epsilon_{i+n}u)^2+(X_{2n+1}^\epsilon X_{2n+1}^\epsilon u)^2)dx.
\end{equation*}
Hence, summing for $l=1$ to $2n+1$ and taking $\tau=\frac{1}{2}$ then the thesis follows.
\end{proof}

The proof above relies crucially on the step-two structure of $\mathbb{H}^n$ and on inequality \eqref{29dicembre}, which in turn depends in an essential way on the precise commutation relations
$[X_i,X_{i+n}]=X_{2n+1}, \qquad i=1,\dots,n.$ In a general step-two Carnot group the commutator structure is more complicated, and establishing analogous mixed estimates would require a separate, more detailed analysis.

The argument in the proof above can be adapted to the case $\beta=0$ in order to obtain the following result.
\begin{cor}\label{2dicembre}
There exists a constant $c=c(q,n,L)>0$, independent of $\epsilon$ and $\delta$, such that, for every weak solution $u\in W_\epsilon^{1,q}(\Omega)$ to \eqref{eqneriemannianapprox} and for every non-negative $\eta\in C_c^\infty(\Omega)$ one has
\begin{equation*}
    \int_\Omega\ \eta^2\ \|\Phi(\nabla_\epsilon u)\|_\lambda^2\ dx\leq c\|\nabla_\epsilon\eta\|_\infty^2\int_{\textnormal{supp}(\eta)} w(\nabla_\epsilon u)^2\ \Lambda dx,
\end{equation*}
where $\|\Phi\|_\lambda^2$ is defined in \eqref{Phi_def} and $w$ and $\Lambda$ are defined in \eqref{26ottobre}.
\end{cor}

\begin{cor}
	Let $2\leq q<\infty$ and $u\in W_\epsilon^{1,q}(\Omega)$ be a weak solution of \eqref{eqneriemannianapprox}, then
	\begin{equation*}
		\lambda_{i,\delta}(\nabla_\epsilon u)X_i^\epsilon u\in HW^{1,2}_{loc}(\Omega),
	\end{equation*}
	for any $i=1,\ldots,2n+1$.
\end{cor}

\begin{proof}
The statement follows from Corollary \ref{2dicembre} and 
\begin{equation*}
	\left|X_j\left(\lambda_{i,\delta}(\nabla_\epsilon u)X_i^\epsilon u\right)\right|^2\leq\|\Phi\|_\lambda^2,
\end{equation*}
for any $i,j=1,\ldots,2n+1$.
\end{proof}

\begin{cor}\label{corollary}
There exists a constant $c=c(q,n,L)>0$, independent of $\epsilon$ and $\delta$, such that, for every weak solution $u\in W_\epsilon^{1,q}(\Omega)$ to \eqref{eqneriemannianapprox}, for every $\alpha\geq\beta\geq2$ and for every non-negative $\eta\in C_c^\infty(\Omega)$ one has
\begin{equation*}
	\int_\Omega\eta^{\alpha+2}\ |X_{2n+1}u|^\beta\ \|\Phi(\nabla_\epsilon u)\|_\lambda^2\ dx\\\leq c^{\frac{\beta}{2}}(\alpha+2)^\beta\|\nabla_\epsilon\eta\|_\infty^\beta\int_\Omega \eta^{\alpha-\beta+2}\ w(\nabla_\epsilon u)^{\beta}\ \|\Phi(\nabla_\epsilon u)\|_\lambda^2 dx.
\end{equation*}
\end{cor}

\begin{proof}	
The statement follows by applying H\"older's inequality with exponents $\frac{\beta}{\beta-2}$ and $\frac{\beta}{2}$ to the right hand-side in Lemma \ref{mainlemma}, and representing $\eta^\alpha= \eta^{\frac{(\alpha+2)(\beta-2)}{\beta}} \eta^{\frac{2 (\alpha-\beta+2)}{\beta}}$. Therefore,
\begin{align*}
	\int_\Omega \eta^{\alpha+2}|X_{2n+1}u|^{\beta}\|\Phi\|_\lambda^2dx\leq
	c(\alpha+2)^2\|\nabla_\epsilon\eta\|_\infty^2\int_\Omega\eta^\frac{(\alpha+2)(\beta-2)}{\beta}\eta^\frac{2(\alpha-\beta+2)}{\beta}w^2|X_{2n+1}u|^{\beta-2}\|\Phi\|_\lambda^2 dx\\\leq c^{\frac{\beta}{2}}(\alpha+2)^\beta\|\nabla_\epsilon\eta\|_\infty^\beta\left(\int_\Omega\eta^{\alpha -\beta+2}w^\beta\|\Phi\|_\lambda^2 dx\right)^{\frac{2}{\beta}}\left(\int_\Omega \eta^{\alpha+2}|X_{2n+1}u|^\beta\|\Phi\|_\lambda^2 dx\right)^{\frac{\beta-2}{\beta}}.
\end{align*}
Dividing both sides by $\left(\int_\Omega \eta^{\alpha+2}|X_{2n+1}u|^\beta\|\Phi\|_\lambda^2 dx\right)^{\frac{\beta-2}{\beta}}$ the thesis follows.
\end{proof}

Now we are in position to get a uniform (in $\delta$ and $\epsilon$) Caccioppoli-type inequality for $\nabla_\epsilon u$, in which the term containing $X_{2n+1}u$ has disappeared.

\begin{teo}\label{maintheorem}
There exists a constant $c=c(q,n,L)>0$, independent of $\epsilon$ and $\delta$, such that, for every weak solution $u\in W_\epsilon^{1,q}(\Omega)$ to \eqref{eqneriemannianapprox}, for every $\alpha\geq0$ and $\beta\geq2$ and for every non-negative $\eta\in C_c^\infty(\Omega)$ one has
\begin{equation*}
	\int_\Omega \eta^{\alpha+2}\ w(\nabla_\epsilon u)^{\beta}\ \|\Phi(\nabla_\epsilon u)\|_\lambda^2\ dx\leq c(\alpha+\beta+2)^{10}K\int_{\Omega} \eta^{\alpha}\ w(\nabla_\epsilon u)^{\beta+2}\ \Lambda dx,
\end{equation*}
where $K=\|\nabla_\epsilon\eta\|^2_\infty+\|\eta X_{2n+1}\eta\|_\infty$, $\|\Phi\|_\lambda^2$ is defined in \eqref{Phi_def} and $w$ and $\Lambda$ are defined in \eqref{26ottobre}.
\end{teo}

\begin{proof}
We apply Lemma \ref{caccioppoliXuriemannian} with $\eta=\eta_1^\frac{\alpha+2}{2} $  
\begin{equation}\label{29gennaio3}
\aligned
	\int_\Omega \eta_1^{\alpha+2}w^\beta\|\Phi\|_\lambda^2dx&\leq c(\beta+1)\left(\frac{\alpha+2}{2}\right)^2K\int_\Omega \eta_1^{\alpha}w^{\beta+2}\Lambda dx\\&+
	c(\beta+1)^2\left(\frac{\alpha+2}{2}\right)^2\int_\Omega \eta_1^{\alpha+2}w^\beta\Lambda|X_{2n+1}u|^2dx\\&\leq c(\alpha+\beta+2)^3K\int_\Omega \eta_1^{\alpha} w^{\beta+2}\Lambda dx+c(\alpha+\beta+2)^4\int_\Omega\eta_1^{\alpha+2}w^\beta\Lambda|X_{2n+1}u|^2dx,
\endaligned
\end{equation}
where $K=\|\nabla_\epsilon\eta_1\|_\infty^2+\|\eta_1X_{2n+1}\eta_1\|_\infty$. 
We apply H\"older's inequality
with exponents $\frac{\beta+2}{\beta}$ and $\frac{\beta+2}{2}$ to the second term, noticing that $\eta_1^{\alpha+2}= \eta_1^\frac{2(\alpha+\beta+2)}{\beta+2} \eta_1^\frac{\alpha\beta}{\beta+2}$. Hence, 
\begin{multline*}
	\int_\Omega\eta_1^{\alpha+2}w^\beta\Lambda|X_{2n+1}u|^2dx\leq \left(\int_\Omega\eta_1^{\alpha + \beta+2}|X_{2n+1}u|^\beta\Lambda|X_{2n+1}u|^2dx\right)^{\frac{2}{\beta+2}}\left(\int_{\Omega}\eta_1^{\alpha}w^{\beta+2}\Lambda dx \right)^{\frac{\beta}{\beta+2}}
\end{multline*}
(we continue by using \eqref{29dicembre})
\begin{equation*}
	\leq\left(\int_\Omega\eta_1^{\alpha+\beta+2}|X_{2n+1}u|^\beta \|\Phi\|_\lambda^2dx\right)^{\frac{2}{\beta+2}}\left(\int_{\Omega}\eta_1^{\alpha}w^{\beta+2}\Lambda dx\right)^{\frac{\beta}{\beta+2}}
\end{equation*}
(we continue by estimating the first term in the right hand side by using Corollary \ref{corollary} with $\alpha_1=\alpha+\beta\ge2$)
\begin{equation*}
	\leq c^{\frac{\beta}{\beta+2}}(\alpha+\beta+2)^{\frac{2\beta}{\beta+2}}\|\nabla_\epsilon \eta_1\|_\infty^{\frac{2\beta}{\beta+2}}\left(\int_\Omega\eta_1^{\alpha+2} w^\beta\|\Phi\|_\lambda^2dx\right)^{\frac{2}{\beta+2}}\left(\int_{\Omega}\eta_1^{\alpha}w^{\beta+2}\Lambda dx\right)^{\frac{\beta}{\beta+2}},
\end{equation*}
where $c=c(q,n,L)>0$ is the constant in Lemma \ref{mainlemma}. Therefore,
calling $$I:=\int_\Omega \eta_1^{\alpha}w^{\beta+2}\Lambda dx,\quad J:=\int_\Omega \eta_1^{\alpha+2}w^\beta\|\Phi\|_\lambda^2 dx$$
and plugging the previous estimate to \eqref{29gennaio3}, we get
\begin{equation*}
	J\leq c(\alpha+\beta+2)^3KI+c(\alpha+\beta+2)^4c^{\frac{\beta}{\beta+2}}(\alpha+\beta+2)^{\frac{2\beta}{\beta+2}}\|\nabla_\epsilon \eta_1\|_\infty^{\frac{2\beta}{\beta+2}}I^{\frac{\beta}{\beta+2}}J^{\frac{2}{\beta+2}}
\end{equation*}
(we continue by using the Young inequality $ab\leq\frac{a}{p}+\frac{b}{p'}$ with exponents $p=\frac{2}{\beta+2}$ and $\frac{\beta}{\beta+2}$)
\begin{equation*}
	\leq c(\alpha+\beta+2)^3KI+\tau c(\alpha+\beta+2)^4\frac{2}{\beta+2}J+\frac{c^2(\alpha+\beta+2)^6\|\nabla_\epsilon \eta_1\|_\infty^2}{\tau^{\frac{2}{\beta}}}\frac{\beta}{\beta+2}I.
\end{equation*}
Choosing $\tau=\frac{1}{2c(\alpha+\beta+2)^4}$ we obtain
\begin{equation*}
	\frac{1}{2}J\leq c(\alpha+\beta+2)^3KI+2^{\frac{2}{\beta}}c^{2+\frac{2}{\beta}}(\alpha+\beta+2)^{6+\frac{8}{\beta}}KI.
\end{equation*}
Therefore, noticing that $6+\frac{8}{\beta}\leq 10$, for any $\beta\geq2$, we have
\begin{equation*}
	J\leq8c^3(\alpha+\beta+2)^{10}KI
\end{equation*}
and the thesis follows.
\end{proof}

\subsection{Uniform $L^\infty$ estimate for the Riemannian gradient of $\epsilon$-weak solutions}

Using the Moser iteration scheme, we obtain a uniform (in $\delta$ and $\epsilon$) bound for the $L^\infty$ norm of $\nabla_\epsilon u.$
 
\begin{teo}[$L^\infty$ estimate for $\nabla_\epsilon u$]\label{linftyestimatesepsilonu}
Let $2\leq q<\infty$ and $u\in W_\epsilon^{1,q}(\Omega)$ be a weak solution of \eqref{eqneriemannianapprox}. Then, for any ball $B_\epsilon (r)$ such that $B_\epsilon (2r)\subset\Omega$ it holds that 
\begin{equation*}
	\|\nabla_\epsilon u\|_{L^\infty(B_\epsilon(r))}\leq c \left(\fint_{B_\epsilon(2r)} \left(\delta+|\nabla_\epsilon u|_\epsilon^2\right)^{\frac{q}{2}}dx\right)^{\frac{1}{q}},
\end{equation*}
where $c=c(q,n,L)>0$.
\end{teo}

\begin{proof}
For any $i=1,\ldots,2n+1$ and any $\beta\geq0$, we compute 
\begin{equation*}
	\nabla_\epsilon\left(\left(\delta+|X^\epsilon_iu|^2\right)^{\frac{q+\beta}{4}}\right)=\left(\frac{q+\beta}{2}\right)\left(\delta+|X^\epsilon_iu|^2\right)^{\frac{q+\beta-4}{4}}X_i^\epsilon u\nabla_\epsilon X^\epsilon_iu.
\end{equation*}
Hence,
\begin{equation}\label{lip1}
    \left|\nabla_\epsilon\left(\left(\delta+|X^\epsilon_iu|^2\right)^{\frac{q+\beta}{4}}\right)\right|_\epsilon^2\leq\frac{1}{4}\left(q+\beta\right)^2\left(\delta+|X_i^\epsilon u|^2\right)^{\frac{q+\beta-2}{2}}|\nabla_\epsilon X^\epsilon_iu|_\epsilon^2.
\end{equation}
	
Since $q\geq2$, from \eqref{lambdaideltaeps} it follows that 
\begin{equation}\label{6gennaio}
    (\delta+|X_i^\epsilon u|^2)^\frac{q-2}{2}\leq(\delta+|X_i^\epsilon u|^2+|X_{i+n}^\epsilon u|^2+|X_{2n+1}^\epsilon u|^2 )^\frac{q-2}{2}= \lambda_{i,\delta}(\nabla_\epsilon u),
\end{equation}
for any $i=1,\ldots,2n+1$. Therefore, \eqref{lip1} and \eqref{6gennaio} imply 
\begin{equation*}
	\left|\nabla_\epsilon\left(\left(\delta+|X^\epsilon_iu|^2\right)^{\frac{q+\beta}{4}}\right)\right|_\epsilon^2\leq\frac{1}{4}\left(q+\beta\right)^2w^{\beta}\|\Phi\|_\lambda^2,\qquad\forall i=1,\ldots,2n+1,
\end{equation*}
where we recall that $w=w(\nabla_\epsilon u)=\left(\delta+|\nabla_\epsilon u|_\epsilon^2\right)^{\frac{1}{2}}$. 
	
Now let $\eta\in C_c^\infty(\Omega)$ be a non-negative cut-off function. Multiplying both sides of the previous inequality by $\eta^2$, integrating over $\Omega$ and summing over $i=1,\ldots,2n+1$, we obtain
\begin{equation*}
	\sum_{i=1}^{2n+1}\int_\Omega \eta^2\left|\nabla_\epsilon\left(\left(\delta+|X^\epsilon_iu|^2\right)^{\frac{q+\beta}{4}}\right)\right|_\epsilon^2dx\leq c_1(q+\beta)^2\int_\Omega \eta^2w^{\beta}\|\Phi\|_\lambda^2dx,
\end{equation*}
where $c_1=c_1(n)>0$. Applying Theorem \ref{maintheorem} with $\alpha = 0$ to the right-hand side, we get
\begin{equation}\label{lip2}
	\sum_{i=1}^{2n+1}\int_\Omega\eta^2 \left|\nabla_\epsilon\left(\left(\delta+|X^\epsilon_iu|^2\right)^{\frac{q+\beta}{4}}\right)\right|_\epsilon^2dx\leq c_2K(\beta+q)^{12}\int_{\text{supp}(\eta)}w^{q+\beta}dx,
\end{equation}
where now $c_2=c_2(q,n,L)>0$ and $K=\|\nabla_\epsilon\eta\|_\infty^2+\|\eta X_{2n+1}\eta\|_\infty$.

Let $B_\epsilon(r)$ be any Riemannian ball. We consider a standard choice of cut-off function $\eta\in C_c^\infty(B_\epsilon(r))$ such that $0\leq \eta\leq 1$ and $\eta\equiv 1$ in $B_\epsilon(r')$, with $0<r'<r$, and
\begin{equation*}
    \|\nabla_\epsilon \eta\|_\infty\leq \frac{4}{r-r'},\qquad \|\nabla^2_\epsilon\eta\|_\infty\leq\frac{16n}{(r-r')^2}.
\end{equation*}
Then
\begin{equation}\label{16dicembre25}
	\sum_{i=1}^{2n+1}\int_{B_\epsilon(r')}\left|\nabla_\epsilon\left(\left(\delta+|X^\epsilon_iu|^2\right)^{\frac{q+\beta}{4}}\right)\right|_\epsilon^2dx\leq \frac{c_3(\beta+q)^{12}}{(r-r')^2}\int_{B_\epsilon(r)}w^{q+\beta}dx,
\end{equation}
where $c_3=c_3(q,n,L)$.

Moreover, a direct computation implies that
\begin{equation*}
	\left(w^{\frac{q+\beta}{2}}\right)^{\frac{2N}{N-2}}\leq c_4\left(2n+1\right)^{\frac{N}{N-2}(q+\beta)}\sum_{i=1}^{2n+1}\left(\left(\delta+|X^\epsilon_iu|^2\right)^{\frac{\beta+q}{4}}\right)^{\frac{2N}{N-2}},
\end{equation*}
with $c_4=c_4(n)>0$. Integrating both sides over $B_\epsilon(r')$, then raising the resulting inequality to the power $\frac{N-2}{N}$, and using the Sobolev inequality (Theorem \ref{sobolevineqriemannian}) with $p=2$ together with \eqref{16dicembre25}, we obtain
\begin{equation}\label{lip3}
\aligned
    \left(\fint_{B_\epsilon(r')} \left(w^{\frac{q+\beta}{2}}\right)^{\frac{2N}{N-2}}dx\right)^{\frac{N-2}{N}}&\leq c_5(n)(2n+1)^{q+\beta}\sum_{i=1}^{2n+1}\left(\fint_{B_\epsilon(r')}\left(\left(\delta+|X^\epsilon_i u|^2\right)^{\frac{\beta+q}{4}}\right)^{\frac{2N}{N-2}}dx\right)^{\frac{N-2}{N}}\\
    &\leq \frac{c_6(2n+1)^{q+\beta}(\beta+q)^{12}}{(r-r')^2}\fint_{B_\epsilon(r)}w^{q+\beta}dx,
\endaligned
\end{equation}
where $c_6=c_6(q,n,L)>0$. If we denote by $k=\frac{N}{N-2}$, we can conclude that
\begin{equation}\label{lip4}
	\left(\fint_{B_\epsilon(r')} w^{(q+\beta)k}dx\right)^{\frac{1}{k}}\leq \frac{c_6(\beta+q)^{12}(2n+1)^{q+\beta}}{(r-r')^2}\fint_{B_\epsilon(r)} w^{q+\beta}dx.
\end{equation}
	
Fix two concentric balls $B_\epsilon(\tau r)\subset B_\epsilon(r)$, with $0<\tau<1$. We consider a sequence of decreasing radii
\begin{equation*}
	r_j=\tau r+\frac{r-\tau r}{2^j}\searrow \tau r,
\end{equation*}
a sequence of non-negative cut-off functions $\eta_j\in C_c^\infty(B_\epsilon(r_j))$ such that $0\leq\eta_j\leq1$, $\eta_j\equiv1$ on $B_\epsilon(r_{j+1})$, and
\[
\|\nabla_\epsilon\eta_j \|_\infty\leq \frac{4}{r_j-r_{j+1}}, \qquad \|\nabla^2_\epsilon\eta_j \|_\infty\leq \frac{16n}{(r_j-r_{j+1})^2},
\]
together with a sequence of increasing exponents
\[
\beta_j=(q+2)k^j-q\geq2.
\]
Using these exponents in \eqref{lip4} and then raising both sides to the power $\frac{1}{(q+2)k^j}=\frac{1}{\beta_j+q}$, we obtain 
\begin{multline}\label{lip5}
	\left(\fint_{B_\epsilon(r_{j+1})} w^{(q+2)k^{j+1}}dx\right)^{\frac{1}{(q+2)k^{j+1}}}\\\leq \left(\frac{c}{(1-\tau)^2}\right)^{\frac{1}{(q+2)k^{j}}}((q+2)k^j)^{\frac{12}{(q+2)k^j}}\left(\fint_{B_\epsilon(r_j)}w^{(q+2)k^j}dx\right)^{\frac{1}{(q+2)k^{j}}},
\end{multline}
where $c=c(q,n,L)>0$. Moreover, if we denote by $\alpha_j=(q+2)k^j$ and iterate \eqref{lip5}, we get 
\begin{equation*}
\aligned
    \left(\fint_{B_\epsilon(r_{m+1})}w^{\alpha_{m+1}}dx\right)^{\frac{1}{\alpha_{m+1}}}&\leq \left(\frac{c}{(1-\tau)^2}\right)^{\sum_{j=0}^m\frac{1}{\alpha_{j}}}\prod_{j=0}^m\alpha_j^{\frac{12}{\alpha_j}}\left(\fint_{B_\epsilon(r)} w^{q+2}dx\right)^{\frac{1}{q+2}}\\
	&\leq \left(\frac{c}{(1-\tau)^2}\right)^{\sum_{j=0}^\infty\frac{1}{\alpha_{j}}}\prod_{j=0}^\infty\alpha_j^{\frac{12}{\alpha_j}}\left(\fint_{B_\epsilon(r)}w^{q+2}dx\right)^{\frac{1}{q+2}}.
\endaligned
\end{equation*}
Now
\begin{equation*}
	\sum_{j=0}^\infty \frac{1}{\alpha_j}=\frac{1}{q+2}\sum_{j=0}^\infty \frac{1}{k^j}=\frac{k}{(k-2)(q+2)}=\frac{N}{2(q+2)},
\end{equation*}
and
\begin{equation*}
	\log\left(\prod_{j=0}^\infty\alpha_j^{\frac{12}{\alpha_j}}\right)=\sum_{j=0}^\infty\frac{12}{\alpha_j}\log\alpha_j=12\frac{\log(q+2)}{q+2}\frac{k}{k-1}+\frac{12}{q+2}\log k\sum_{j=0}^\infty \frac{j}{k^j},
\end{equation*}
which depend only on $q$ and $n$. Hence
\begin{equation*}
	\left(\fint_{B_\epsilon(r_{m+1})} w^{\alpha_{m+1}}dx\right)^{\frac{1}{\alpha_{m+1}}}\leq \left(\frac{c}{(1-\tau)^{\frac{N}{q+2}}}\right)\left(\fint_{B_\epsilon(r)} w^{q+2}dx\right)^{\frac{1}{q+2}}.
\end{equation*}
Since $\alpha_m\to\infty$ as $m\to\infty$, and the averages on the left-hand side tend to the essential supremum of the integrand, we get
\begin{equation*}
	\sup_{B_\epsilon(\tau r)}w\leq\left(\frac{c}{(1-\tau)^{\frac{N}{q+2}}}\right)\left(\fint_{B_\epsilon(r)} w^{q+2}dx\right)^{\frac{1}{q+2}},
\end{equation*}
where $c=c(q,n,L)>0$, and this holds for all $B_\epsilon(r)\subset\Omega$ and for all $0<\delta<1$. 
	
Another iteration argument (see, for instance, \cite[Theorem 5.1]{ricciotti} or \cite[Lemma 3.38]{heinonen2018nonlinear}) implies that
\begin{equation}\label{26gennaio2}
	\sup_{B_\epsilon(\tau r) }w\leq\left(\frac{c'}{(1-\tau)^{\frac{N}{s}}}\right)\left(\fint_{B_\epsilon(r)} w^sdx\right)^{\frac{1}{s}},
\end{equation}
for any $s>0$, where $c'=c'(s,q,n,L)>0$. The thesis follows by choosing $s=q$.
\end{proof}

\subsection{Local Lipschitz regularity for solutions}

Passing to the limit (in $\delta$ and $\epsilon$), we obtain a bound on the $L^\infty$ norm for the horizontal gradient of solutions to \eqref{eq}.

\begin{teo}[$L^\infty$ estimate for the gradient of solutions]

Let $2\leq q<\infty$ and let $u\in HW^{1,q}(\Omega)$ be a weak solution of \eqref{eq}, then $\nabla_Hu\in L^\infty_{loc}(\Omega)$. Moreover, for any ball $B_{SR}(r)$ such that $B_{SR}(2r)\subset\Omega$ it holds 
\begin{equation*}
	\|\nabla_Hu\|_{L^\infty(B_{SR}(r))}\leq c \left(\fint_{B_{SR}(2r)} |\nabla_H u|_H^{q}dx\right)^{\frac{1}{q}},
\end{equation*}
where $c=c(q,n,L)>0$.
\end{teo}

\begin{proof}
We fix a ball $B_{SR}(x_0,r)\subset\Omega$. We may assume that $x_0=0$, and we write $B_{SR}(r)=B_{SR}(0,r)$. We consider a sequence $\left(\psi_\epsilon\right)_{\epsilon>0}\subseteq C^{\infty}(B_{SR}(r))$, such that $\psi_\epsilon\to u$ in $HW^{1,q}(B_{SR}(r))$ for $\epsilon\to0$. Its existence follows from the density of $C^\infty(B_{SR}(r))\cap HW^{1,q}(B_{SR}(r))$ in $HW^{1,q}(B_{SR}(r))$, see \cite[Theorem 1.2.3]{franchi1997approximation}. Standard PDEs arguments ensure the existence and uniqueness of the solution to the following Dirichlet problem
\begin{equation}\label{27gennaio}
\begin{cases}
    \sum_{i=1}^{2n+1}X_i^\epsilon({\Dfdelta i}(\nabla_\epsilon v))=0 \quad&\text{in }B_{SR}(r),\\
	v=\psi_\epsilon\quad &\text{in }\partial B_{SR}(r),
\end{cases}
\end{equation}
see for instance \cite{heinonen2018nonlinear}.
    
Let $u_\epsilon$ be such a solution. By virtue of classical elliptic theory, see for instance \cite{Uraltsevabook}, it follows that $u_\epsilon\in C^\infty(B_{SR}(r))$. By definition of weak solution it follows that 
\begin{equation}\label{weaksol}
    \sum_{i=1}^{2n+1}\int_{B_{SR}(r)}{\Dfdelta i}(\nabla_\epsilon u_\epsilon))X_i^\epsilon (u_\epsilon-\psi_\epsilon)dx=0,
\end{equation}
then \eqref{structurecondition02} implies that 
\begin{equation}
    \int_{B_{SR}(r)}|\nabla_\epsilon u_\epsilon|_\epsilon^q dx\leq \sum_{i=1}^{2n+1}\int_{B_{SR}(r)}{\Dfdelta i}(\nabla_\epsilon u_\epsilon)X_i^\epsilon u_\epsilon dx\leq\int_{B_{SR}(r)}|\nabla_\epsilon u_\epsilon|_\epsilon^{q-1}|\nabla_\epsilon\psi_\epsilon|_\epsilon dx. 
\end{equation}
Hence, by H\"older inequality it follows that there exists $M_0>0$, independent of and $\epsilon$, such that
\begin{equation}\label{17maggio2024}
    \|\nabla_\epsilon u_\epsilon\|_{L^q(B_{SR}(r))}\leq M_0.
\end{equation} 
In particular, $\|u_\epsilon\|_{W_\epsilon^{1,q}(B_{SR}(r))}$ is uniformly bounded with respect to $\epsilon$. Therefore, up to subsequences, $u_\epsilon\rightharpoonup u_0$ in $W_\epsilon^{1,q}(B_{SR}(r))$.
Moreover, since the Riemannian balls $B_\epsilon(r)$ converge in the Hausdorff metric to $B_{SR}$, then there exists $\tau\in(0,1)$ such that $B_\epsilon (2\tau r)\subset B_{SR}(r)$, for $\epsilon$ small enough. We can apply Theorem \ref{linftyestimatesepsilonu} to infer 
\begin{equation}\label{25gennaio}
   \|\nabla_\epsilon u_\epsilon\|_{L^\infty(B_\epsilon(\tau r))}\leq c \left(\fint_{B_\epsilon(2\tau r)} \left(\delta+|\nabla_\epsilon u_\epsilon|_\epsilon^2\right)^{\frac{q}{2}}dx\right)^{\frac{1}{q}},
\end{equation}
where $c=c(q,n,L)>0$. 

Moreover, from Theorem \ref{maintheorem} and \eqref{17maggio2024} it follows that there exists $M_1>0$, independent of $\epsilon$ and $\delta$, such that 
$$\|\nabla_\epsilon\left(|X_i^\epsilon u_\epsilon|^{q+\beta}\right)\|_{L^2(B_\epsilon(\tau r))}\leq M_1,\quad \forall \beta\geq2.$$
The Rellich–Kondrachov theorem implies that $|X_i^\epsilon u_\epsilon|^{q+\beta}$ converges, up to subsequences, in $L^\sigma(B_\epsilon(\tau r))$, for any $\sigma<2^*$. Hence, $X_i^\epsilon u_\epsilon$ converges, up to subsequences, in $L^s(B_\epsilon(\tau r))$ for any $s\in(1,\infty)$, from which it follows that $X_i^\epsilon u_\epsilon$ converges, up to subsequences, point-wise a.e. in $B_\epsilon(\tau r)$. Hence, by definition of weak derivative, we can conclude that $\nabla_\epsilon u_\epsilon\to\nabla_Hu_0$ a.e. in $B_\epsilon(\tau r)$, up to subsequences. 

On the one hand, we can pass to the limit in \eqref{25gennaio}, both in $\epsilon$ and $\delta$, and get
\begin{equation*}
   \|\nabla_H u_0\|_{L^\infty(B_{SR}(\tau r))}\leq c \left(\fint_{B_{SR}(r)} |\nabla_H u_0|_H^qdx\right)^{\frac{1}{q}};
\end{equation*}
on the other hand we can pass to the limit in the weak formulation of \eqref{27gennaio} and we get that $u_0$ is a weak solution to \eqref{eq}. Since $u_\epsilon-\psi_\epsilon\in HW^{1,q}_0(B_{SR}(r))$, and this space is closed under weak convergence, we get that $u_0-u\in HW^{1,q}_0(B_{SR}(r))$. Then, the comparison principle implies that $u_0=u$ in $B_{SR}(r)$. This estimates holds for all $B_{SR}(r)\subset\Omega$. Now the thesis follows from a simple covering argument.
\end{proof}

Therefore, the following corollary holds.

\begin{cor}[Local Lipschitzianity]
Let $2\leq q<\infty$ and $u\in HW^{1,q}(\Omega)$ be a weak solution of \eqref{eq}, then $u$ is locally Lipschitz continuous in $\Omega$. Moreover, for any ball $B_{SR}(r)$ such that $B_{SR}(2r)\subset\Omega$ it holds that 
\begin{equation*}
	|u(x)-u(y)|\leq c \left(\fint_{B_{SR}(2r)} |\nabla_H u|_H^{q}dz\right)^{\frac{1}{q}}d_{SR}(x,y),
\end{equation*}
where $c=c(q,n,L)>0$.
\end{cor}

\section*{Acknowledgements}
M. C. and G.C. are supported by the project PRIN 2022 F4F2LH - CUP J53D23003760006, G.C. is funded by MNESYS PE12 (PE0000006).

A.C. acknowledges partial support from grants PID2020-112881GB-I00, PID2021-125021NAI00 (Spanish Government) and 2021-SGR-00071 (Catalan Government).

\printbibliography

\end{document}